\documentclass[11pt]{amsart}
\usepackage{amsfonts}
\usepackage{amsmath}
\usepackage{amssymb}
\usepackage{amsthm}
\usepackage{ color, ulem}
\usepackage{mathtools}

\usepackage[dvipsnames]{xcolor}
\usepackage{soul}

\input xy
\xyoption {all}

\newcommand{\comment}[1]{}
\newtheorem{theorem}{Theorem}
\newtheorem {lemma}{Lemma}

\newtheorem {proposition}{Proposition}
\newtheorem{conjecture}{Conjecture}

\theoremstyle{definition}
 
\theoremstyle {definition}
\newtheorem{remark}{Remark}
\newtheorem{definition}{Definition}

\begin{document}
\baselineskip=16.25pt
\begin{center}
\title[Higher rank Segre integrals]{Higher rank Segre integrals over the Hilbert scheme of points}
\author{A. Marian}
\address{Department of Mathematics, Northeastern University} 
\email {a.marian@neu.edu}
\author{D. Oprea}
\address{Department of Mathematics, University of California, San Diego}
\email {doprea@math.ucsd.edu}
\author{R. Pandharipande}
\address{Department of Mathematics, ETH Z\"urich}
\email {rahul@math.ethz.ch}
\maketitle

\end{center}
\begin{abstract}{Let $S$ be a nonsingular projective surface. Each vector bundle $V$ on $S$ of rank $s$ induces a tautological vector bundle over the Hilbert scheme of $n$ points of $S$. When $s=1$, the top Segre classes of the
tautological bundles are given by a recently proven formula conjectured in 1999 by M. Lehn. We calculate here the Segre classes of 
the tautological bundles for all ranks $s$ over all $K$-trivial surfaces. Furthermore, in rank $s=2$, the Segre integrals are determined for all surfaces, thus establishing a full analogue of Lehn's formula. We also give conjectural formulas for certain series of Verlinde Euler characteristics over the Hilbert schemes of points.}
\end{abstract}

\vskip.1in

\section{Introduction}

\subsection {The Hilbert scheme of points} Let $S$ be a nonsingular projective surface, and let $S^{[n]}$ denote the Hilbert scheme of $n$ points on $S$. Each line bundle $L\to S$ gives rise to a tautological rank $n$ vector bundle $L^{[n]}\to S^{[n]}$ via the assignment $$\zeta\mapsto H^0(L\otimes \mathcal O_\zeta).$$ 

Tautological integrals over the Hilbert scheme of points or over their geometric subsets emerge often in enumerative geometry. We mention three situations studied previously:
\begin {itemize}
\item [(i)] the count of $n$-nodal curves in fixed linear system $\mathbb P^n\subset |L|$ where $L\to S$ is sufficiently positive. These counts were recast by G\"{o}ttsche \cite {G1} in terms of Hilbert schemes. Specifically, for a suitable geometric subscheme $$W_n\hookrightarrow S^{[3n]}$$ birational to $S^{[n]}$, the Severi degrees are encoded by the series $$\mathsf G(z)=\sum_{n=0}^{\infty} z^n \int_{W_n} c(L^{[n]}).$$ 
\item [(ii)] the Verlinde series of holomorphic Euler characteristics \cite{EGL}: $$\mathsf V(z)=\sum_{n=0}^{\infty} z^n \cdot \chi(S^{[n]}, L_{(n)}\otimes E^r).$$ Here the line bundle $()_{(n)}$ is pulled back from the symmetric product, and $E$ denotes $-\frac{1}{2}$ of the exceptional divisor on the Hilbert scheme. 
\item [(iii)] integrals of Segre classes of tautological bundles considered by Lehn \cite {L}: $$\mathsf S(z)=\sum_{n=0}^{\infty} z^n \int_{S^{[n]}} s(L^{[n]}).$$ 
\end {itemize}

There are strong reasons to study all three series above. The count of curves in (i) is a central enumerative question. The Verlinde series (ii) captures all holomorphic Euler characteristics of line bundles over the Hilbert scheme of points. Calculating the series has applications to the study of $K$-theoretic Donaldson invariants on moduli spaces of sheaves of arbitrary rank, as pursued in \cite{gny}, \cite{gy}, \cite{g3}. Indeed over elliptic surfaces, higher-rank invariants can in many cases be matched to Verlinde numbers on the Hilbert scheme via Fourier-Mukai techniques. Furthermore, knowledge of the series (ii) is crucial to the strange duality problem for sheaves on surfaces. Interesting variants of the Verlinde series, involving symmetric and exterior powers of tautological bundles, were recently studied in \cite{a}; see also \cite {Sc}.

The main focus of this paper is the Segre series (iii) in its general form, when $L$ is replaced by an arbitrary higher rank vector bundle. The Segre integrals appeared first in the study of rank $2$ Donaldson invariants of rational surfaces \cite {T, Ty}. 
More recently, the Segre series have turned up in the intersection theory of various parameter spaces of sheaves with sections (higher rank Quot schemes, higher rank stable pairs), 
where they naturally enter localization calculations via the obstruction bundles (cf \cite{mop}). This is a setup which parallels the appearance of Hodge integrals in Gromov-Witten theory. The importance of systematically studying the Segre series (iii) and its higher rank generalizations was highlighted in \cite{OP} in the case of curves.

In a significant development, a conjectural relationship between the Verlinde and Segre series was proposed in \cite{J}, aligned with the vaster conjectural framework of strange duality  for moduli spaces of sheaves. This picture acquired further precision in \cite{mop2} where the series $\mathsf S$ and $\mathsf V$ were conjecturally matched through an explicit change of variables.  As a consequence, the Segre series holds the key to understanding the Verlinde numbers (ii). This point will be addressed in Section \ref{vgf}. 

The common feature of the expressions $\mathsf G$, $\mathsf V,$ $\mathsf S$ is that all three factor (cf. \cite {G1}, \cite {Tz}, \cite {KST}, \cite {EGL}) as products of four universal series $${\mathsf U}_1 (z)^{L^2}  \cdot {\mathsf U}_2(z)^{\chi(\mathcal O_S)} \cdot {\mathsf U}_3 (z)^{L\cdot K_S} \cdot {\mathsf U}_4(z)^{K_S^2},$$ with $$\mathsf U_1, \,\mathsf U_2, \mathsf U_3, \mathsf U_4\in \mathbb Q[[z]].$$ The series $\mathsf U_1$ and $\mathsf U_2$ are uniquely determined by $K$-trivial geometries and are typically more accessible. They are known in examples (i)-(iii), see \cite {G1}, \cite {EGL}, \cite {mop}.
For all types of tautological integrals considered over the past few decades, closed formulas for the remaining series 
$\mathsf U_3,\mathsf U_4$ have proven very difficult to calculate or even conjecture. Remarkably, the Segre geometry (iii) is the only nontrivial case known so far; the details are explained in Section \ref{sect13}. 

\subsection{Higher rank} We undertake the study of the Segre series in higher rank. To set the stage, let $V\to S$ be a vector bundle of rank $s$, inducing a tautological rank $sn$ vector bundle $$V^{[n]}\to S^{[n]}$$ over the Hilbert scheme via the assignment $$\zeta\mapsto H^0(V\otimes \mathcal O_\zeta).$$ The construction extends to $K$-theory. For each $K$-theory class $\alpha\in K(S)$, there is an associated $K$-theory class $$\alpha^{[n]}\to S^{[n]}$$ defined via locally free resolutions. We are interested in the Chern numbers of $\alpha^{[n]}$, but it will be more convenient to work with Segre classes. Fixing $\alpha\in K(S)$, we consider its associated Segre series 
\begin{eqnarray}\label{formm2}\mathsf {S}_{\alpha} (z) &=& \sum_{n=0}^{\infty} z^n \int_{S^{[n]}} s(\alpha^{[n]})\\ & =& A_0(z)^{c_2(\alpha)} \cdot A_1 (z)^{c_1(\alpha)^2}  \cdot A_2(z)^{\chi(\mathcal O_S)} \cdot A_3 (z)^{c_1(\alpha)\cdot K_S} \cdot A_4(z)^{K_S^2}.\nonumber \end{eqnarray} The above factorization follows by \cite {EGL}. The five series $$A_0,\, A_1,\, A_2,\, A_3,\, A_4\in \mathbb Q[[z]]$$ are independent of the surface $S$, and depend on $\alpha$ only through the rank $s$. 

\subsection {Lehn's conjecture} \label{sect13}Let us assume first that $\text{rank }\alpha=1$ and $\alpha$ is represented by a line bundle $$L\to S.$$ Then, the series $A_0$ is absent from \eqref{formm2}. Lehn \cite{L} 
 conjectured closed formulas for {\it all} four series $A_1, A_2, A_3, A_4$. In the rephrasing of \cite {mop}, after the change of variables $$z=t(1+2t)^2,$$ the universal Segre series are given by 
\begin{eqnarray*}
A_1(z)&=&(1+2t)^{\frac{1}{2}}\, ,\\ 
A_2(z)&=&(1+2t)^{\frac{3}{2}}\cdot (1+6t)^{-\frac{1}{2}}\, , \\
A_3(z)&=&\frac{1}{2}\cdot (1+2t)^{-1}\cdot \left(\sqrt{1+2t}+\sqrt{1+6t}\right)\,, \\
 A_4(z)&=&4\cdot (1+2t)^{\frac{1}{2}}\cdot (1+6t)^{\frac{1}{2}}\cdot \left(\sqrt{1+2t}+\sqrt{1+6t}\right)^{-2}\, . 
 \end{eqnarray*}
 
The series $A_1$ and $A_2$ were confirmed via the study of the virtual geometry of the Quot scheme of $K$-trivial surfaces in \cite {mop}. Using Reider type techniques, 
Voisin \cite {voisin} proved the vanishing of certain Segre integrals in the case of the blowup of a $K3$ surface in one point, and showed these vanishings determined uniquely the series $A_3, A_4$. 
The series were shown to have Lehn's conjectured expressions via the residue calculations of \cite {mop2}. As a result, the Segre series are calculated for all line bundles over nonsingular projective surfaces $S$.

\subsection {The series $A_0, A_1, A_2$ in arbitrary rank} 
Our first result gives a simple closed form expression for the Segre series of all $K$-trivial surfaces, in all ranks. We prove

\begin{theorem}\label{t1} Let $S$ be a $K$-trivial surface, let $\alpha$ be a $K$-theory class of rank $s$, and set $r=s+1$. We have 
\begin{equation}\label{eq1}\mathsf {S}_{\alpha} (z) =  \sum_{n=0}^{\infty} z^n \int_{S^{[n]}}  s (\alpha^{[n]}) = A_0(z)^{c_2(\alpha)}\cdot A_1(z)^{c_1^2(\alpha)}\cdot A_2(z)^{\chi(\mathcal O_S)}\end{equation}
where the formulas for the series $A_1$, $A_2$, and $A_3$ are  
\begin{eqnarray*}
A_0(z) &=& (1+rt)^{-r} \cdot (1 + (1+r) t)^{r - 1}\, , \\ 
A_1 (z) & = & (1+rt)^{\frac{r-1}{2}} \cdot (1 + (1+r) t)^{-\frac{r}{2}+1}\, ,\\
A_2 (z) & = & (1+rt)^{\frac{r^2-1}{2}} \cdot (1 + (1+r)t)^{-\frac{r^2}{2}+r} \cdot (1 + r(1+r)t)^{-\frac{1}{2}},
\end{eqnarray*} 
after the change of variables
\begin{equation}
z = t \,(1+rt)^r\, .
\label{changez}
\end{equation}

\end{theorem}

\vskip.1in

\begin{remark} 
In \cite {J}, Johnson computed the series $A_0, A_1, A_2$ up to order $6$ and conjectured several connections between them. Building on these calculations, the above formulas for  $A_0, A_1, A_2$ were 
proposed in \cite{mop2}. 
\end{remark} 

\begin{remark} The simplicity of the expressions in Theorem \ref{t1} is deceiving. The individual Segre integrals are quite complicated. For instance, for a $K3$ surface $S$ and $n=2$, we have 
\begin{eqnarray}\label{2pt}\int_{S^{[2]}} s_4(V^{[2]})&=&\frac{1}{4}\left(\left(2 (c_1^2)^2+2c_2^2 -4c_1^2c_2 - 8 c_1^2 +  6 c_2\right)  +s(- 9 c_1^2 + 
 6 c_2 +12)\right.{}\\ \nonumber& &+\left.{}s^2(- 3 c_1^2  + 2 c_2 +22)+12 s^3 + 2 s^4\right)\, ,\end{eqnarray} 
where $V$ is a vector bundle of rank $s$ with Chern classes $c_1$ and $c_2$.
For arbitrary surfaces, the answers are more involved. \end{remark}

\begin{remark} Even in rank $1$, Theorem \ref{t1} gives additional information not covered by Lehn's conjecture. While the series $A_1, A_2$ were predicted by Lehn's formula, $$A_0(z)=(1+2t)^{-2}\cdot (1+3t) \text{ for } z=t(1+2t)^2$$ was absent from the line-bundle setup of the conjecture. \end{remark}

\begin{remark} For nonsingular projective curves, the generating series of higher rank Segre integrals was determined in \cite {mop2}. A change of variables similar to \eqref{changez} was needed to express the answer in closed form.  
\end{remark} 

\begin{remark} The expressions $A_1, A_2$ determined by Theorem \ref{t1} were subsequently connected by Mellit with certain generating series of Hurwitz numbers in genus zero \cite {M}.

\end{remark}

\subsection {The series $A_3, A_4$.} For surfaces $S$ which are not $K$-trivial, and an arbitrary class $\alpha \in K (S)$, the universal series $A_3, A_4$ appear in the Segre generating
function $\mathsf {S}_{\alpha} (z)$ and are very difficult to evaluate. \vskip.1in

We calculate the two remaining universal series for $\text{rank} \, \alpha = 2$, thus proving a full 
rank-two analogue of Lehn's formulas. Thus, we provide a second nontrivial situation for which all series involved are determined. Using the change of variables $$z=t(1+3t)^3,$$ Theorem \ref{t1} gives
\begin{eqnarray*}
A_0(z) &=& \frac{(1 + 4 t)^{2}}{ (1+3t)^{3}}\, , \\ 
A_1 (z) & = & \frac{(1+3t)}{(1 + 4 t)^{\frac{1}{2}}}\, ,\\
A_2 (z) & = & \frac{(1+3t)^{4}}{ (1 + 4t)^{\frac{3}{2}} \cdot (1 + 12t)^{\frac{1}{2}}}.
\end{eqnarray*}

To describe the remaining two series, let $\mathsf y(t)$ be the unique real solution of the quartic equation $$\frac{y\cdot(1 + y)^2}{(1 - y)(1-y^3)}= \frac{t}{1+3t},\,\,\, y(0)=0.$$ 
       We have $$\mathsf y(t)=t - 6 t^2 + 41 t^3 - 314 t^4 + 2630 t^5+\ldots.$$ 
              \begin{theorem} \label{conj2} For $\text{rank }\alpha=2$ the remaining universal Segre series are \begin{eqnarray*}A_3(z) &=& (1+3t)^{-1} \cdot \left(\frac{\mathsf y(t)}{t}\right)^{-1/2}\\
A_4(z)&=&(1+3t)\cdot \left(\frac{\mathsf y(t)}{t}\right)^{3}\cdot \frac{(1+\mathsf y)^2}{1-\mathsf y}\cdot \frac{1}{\mathsf y'}.\end{eqnarray*}\end{theorem}

Based on these expressions, we also formulate 

\begin{conjecture}\label{co7} For all ranks, 
the functions $A_3, A_4$ are algebraic. 
\end{conjecture} 

\begin{remark}
Finding the series $A_3$ and $A_4$ for $\alpha$ of arbitrary rank remains an open question. In addition to Lehn's formula for rank one, and the substantial rank two analysis carried out here, 
a few cases are trivially known:
\begin{itemize} 
\item the case $\text{rank } \alpha=-1$ is immediate. Set $\alpha=-L$ for some line bundle $L\to S$. For dimension reasons we have $$\int_{S^{[n]}} s(\alpha^{[n]})=\int_{S^{[n]}} c((-\alpha)^{[n]})=0 \text{ for } n\geq 1\implies A_3(z)=A_4(z)=1.$$ 
\item for $\text{rank }\alpha=-2$, we can prove geometrically, see \eqref {rank2222} for instance, that $$\int_{S^{[n]}} s(\alpha^{[n]})=\int_{S^{[n]}} c((-\alpha)^{[n]})=\binom{c_2(-\alpha)}{n}.$$ We again obtain $$A_3(z)=A_4(z)=1.$$ 
\item for $\text{rank }\alpha=0$, we have $A_4(z)=1$ just by setting $\alpha=0$. Regarding the remaining series, we state the following 
\end{itemize}
\begin{conjecture} \label{c3}In rank $0$, we have $$A_3(z)=(1+t)^{-1}\cdot (1+2t)^{\frac{1}{2}} \text { for } z=t(1+t).$$\end{conjecture} 

\end{remark}

\vskip.1in

\subsection {The Verlinde generating function} \label{vgf} We now turn to the generating series of holomorphic Euler characteristics of tautological line bundles over the Hilbert scheme. We set $$\mathsf V_{\alpha} (w)=\sum_{n=0}^\infty \chi(S^{[n]}, (\det \alpha)_{(n)}\otimes E^r)\cdot w^n\, ,$$ where the line bundle $()_{(n)}$ is pulled back from the symmetric product, and $E$ denotes $-\frac{1}{2}$ of the exceptional divisor on the Hilbert scheme. Since $\mathsf V_\alpha$ only depends on $\det \alpha$, only four power series are needed in the factorization $$\mathsf V_{\alpha}(w)= B_1(w)^{\chi(c_1(\alpha))}\cdot B_2(w)^{\chi(\mathcal O_S)} \cdot B_3(w)^{c_1(\alpha)\cdot K_S-\frac{1}{2}K_S^2}\cdot B_4(w)^{K_S^2}.$$ The different form of the exponents used here relative to \eqref {formm2} is justified by the fact that with the current choices, we have the following symmetries \begin{equation} \label{symmetry0} B_3 \mapsto B_3^{-1} \text{ as } r\to -r\end{equation} $$B_4 \text{ is invariant under the symmetry } r\to -r.$$ This was noted in Theorem 5.3 in \cite {EGL} as a consequence of Serre duality. Also by \cite {EGL}, we explicitly know $B_1$ and $B_2$. For $$w=t(1+t)^{r^2-1},$$ we have $$B_1(w)=1+t, \,\,\,B_2(w)=\frac{(1+t)^{\frac{r^2}{2}}}{(1+r^2t)^{\frac{1}{2}}}.$$ 

The series $B_3, B_4$ remain mysterious in general. The cases $r = 0$ and $r = \pm 1$ are immediate exceptions, obtained in \cite {EGL}, Lemma 5.1.
Indeed for $r=0$, we have $$\mathsf V_{\alpha}(w)=(1-w)^{-\chi(c_1(\alpha))}\implies B_3(w)=B_4(w)=1,$$ and for $r=\pm 1$, $$\mathsf V_{\alpha}(w)=(1+w)^{\chi(c_1(\alpha))}\implies B_3(w)=B_4(w)=1.$$

For all surfaces, motivated by strange duality, Johnson \cite{J}
predicted a connection between the Segre and Verlinde series $$\mathsf S_{-\alpha} (z)\leftrightarrow \mathsf V_{\alpha}(w),$$ under an undetermined change of variables $z\leftrightarrow w$, and the shift \begin{equation}\label{shiftrank}r=\text{rank }\alpha- 1.\end{equation} These predictions were made precise in \cite {mop2} where the unknown change of variables was proposed: $$z=t(1-rt)^{-r},\,\,\,\, w=\frac{t(1-(r-1)t)^{r^2-1}}{(1-rt)^{r^2}}.$$ 

Using Theorem \ref{t1} we show under suitable numerics an equality between Verlinde numbers and tautological Chern integrals on Enriques surfaces, see Proposition \ref{ccqq}. This is consistent with strange duality.

Since the Segre series is now calculated for rank $\alpha = 1$ (Lehn's formula) and rank $\alpha = 2$ (in Theorem \ref{conj2}), this precise conjectural Verlinde-Segre relationship gives further predictions for the unknown Verlinde universal series $B_3$ and $B_4$. 

For $r=\pm 2$, the formulas for $B_3$ and $B_4$ are captured by Conjecture $2$ of \cite {mop2}. 
Indeed for $r= 2,$ setting \footnote{To be in agreement with the formulas of \cite {mop2}, one has to change $t\mapsto t/(2+2t).$} $$w=t(1+t)^3,$$ we expect  \begin{eqnarray*}B_3(w)&=& \frac{1+\sqrt{1+4t}}{2(1+t)}\\ B_4(w)&=&(1+t)^{\frac{1}{2}}\cdot (1+4t)^{\frac{1}{2}}\cdot \left(\frac{1+\sqrt{1+4t}}{2}\right)^{-\frac{5}{2}}.\end{eqnarray*} 

Thanks to Theorem \ref{conj2}, we obtain the following predictions about the Verlinde series for $r=\pm 3$. Let $\mathsf Y$ be the unique real solution of the quartic equation $$\frac{y\cdot(1 + y)^2}{(1 - y)(1-y^3)}=t,\,\,\,\, y(0)=0.$$ We have $$\mathsf Y(t)=\mathsf y\left(\frac{t}{1-3t}\right)=t-3t^2+14t^4-80t^4+509t^5-3459t^6+\ldots.$$ 

\begin{conjecture} For the Verlinde series $\mathsf V_{\alpha}$ with $r=3$, setting $$w={t(1+t)^8}$$ we have \begin{eqnarray*}B_3(w)&=&(1+t)^{-\frac{3}{2}} \cdot \left(\frac{\mathsf Y(t)}{t}\right)^{-\frac{1}{2}},\\ B_4(w)&=&(1+t)^{\frac{3}{4}}\cdot \left(\frac{\mathsf Y(t)}{t}\right)^{\frac{13}{4}}\cdot \frac{(1+\mathsf Y)^2}{1-\mathsf Y}\cdot \frac{1}{\mathsf Y'}.\end{eqnarray*} The expressions for $r=-3$ are obtained via the symmetry \eqref{symmetry0}.
\end {conjecture}
Finding the general expression for the unknown series $B_3$ and $B_4$ for arbitrary $r$, thus determining all rank $1$ Verlinde numbers, is a central question in the enumerative theory of Hilbert schemes of points on surfaces.\vskip.1in

\begin{remark}

Due to the rank shift \eqref{shiftrank}, the Serre duality symmetry $$s\to -s$$ on the Verlinde side 
translates into a {\it conjectural} transformation rule for the remaining unknown Segre universal series $A_3, A_4,$ as $$s\to -s-2.$$ In particular, we can also state predictions for the Segre series $A_3, A_4$ in $\text{rank }\alpha=-3$ and $\text{rank } \alpha=-4$ from the series in $\text{rank } \alpha=1$ and $\text{rank }\alpha=2$. The exact expressions are more cumbersome, but they support Conjecture \ref{co7}. 
\end{remark} 

\subsection {Strategy of proofs.} 
\label{sttg}
We first explain how Theorem \ref{t1} is derived. To find the three universal series $A_0, A_1, A_2$, we pick $\alpha$ to be the class of a suitable vector bundle $$V\to S$$ over a $K3$ surface $S$. As witnessed by \eqref{2pt}, the Segre integrals are generally very complicated. It is a key observation that the answers take a simpler form for sheaves with small deformation spaces. The most beautiful formulas are obtained for spherical vector bundles $V$ and for vector bundles with isotropic Mukai vectors. Let $$v = \text{ch} \, V \cdot \sqrt{\text{td} (S)} \in H^{2\star} (S, \mathbb Z)$$ be the Mukai vector of $V$, and recall the Mukai pairing $$\langle v, v\rangle=\int_{S} v_2^2-2v_0v_4 \quad \text{ where } v=(v_0, v_2, v_4)\in H^{2\star}(S).$$
Set $\chi=\chi(S, V).$ We show:
\begin{theorem} Let $S$ be a K3 surface, and let $V\to S$ be a
 rank $s=r-1$ vector bundle.
\begin{itemize}\label{p1}
\item [(i)] If $\langle v, v\rangle = -2$, then $$\int_{S^{[n]}}s_{2n}(V^{[n]})= r^n \binom{\chi - rn}{n}\, .$$
\item [(ii)] If $\langle v, v\rangle = 0$, then $$\int_{S^{[n]}}s_{2n}(V^{[n]})=r^n \left ( -r+\frac{1}{r} +\frac{\chi}{n} \right ) \binom{\chi-rn-1}{n-1}\, .$$
\end {itemize}
\end {theorem} 

Theorem \ref{p1} is proven by a direct geometric argument using Reider techniques \cite{R}, also employed in rank $1$ in \cite{voisin}. The crucial insight here is the identification of the optimal geometric setup for which the complicated Segre integrals become manageable. We will then show by a residue calculation that Theorems \ref{t1} and \ref{p1} are  equivalent.
\vskip.1in

The analysis is quite intricate for Theorem \ref{conj2}, requiring in particular delicate excess calculations for Segre classes. The key statements are captured by the following
\begin{theorem}\label{t4} Let $\pi: S\to X$ be the blowup of a $K3$ surface $X$ at a point, with exceptional divisor $E$. Let $V_0\to X$ be a rank $2$ bundle whose Mukai vector satisfies $\langle v_0, v_0\rangle=-2$. Set $$V=\pi^{\star} V_0\otimes E^{-k},$$ and assume that $$\chi(V)=4n-1.$$
\begin{itemize}
\item [(i)] If $k=n-2$, then $$\int_{S^{[n]}} s_{2n}(V^{[n]})=(-1)^n(2n+1).$$
\item [(ii)] If $k=n-1$, we have $$\int_{S^{[n]}} s_{2n}(V^{[n]})=1 \, \, \, \text{for} \, \, \, n\equiv 0\mod 3.$$
$$\int_{S^{[n]}} s_{2n}(V^{[n]})=0 \, \, \, \text{for} \, \, \, n\not \equiv 0 \mod 3.$$
\end {itemize}
\end{theorem}
We finally show by a residue calculation that Theorems \ref{conj2} and \ref{t4} are equivalent. 

{\subsection {Moduli of surfaces} The Segre integrals can be viewed as part of a richer theory which is important to explore further. For each flat family
$$\pi:S\to B\, $$
of nonsingular projective surfaces carrying line bundles
$L_1,\ldots, L_\ell\to S$, we define the $\kappa$-classes $$\kappa[a_1, \ldots, a_\ell, b]=\pi_{\star} \left(c_1(L_1)^{a_1}\cdots c_1(L_\ell)^{a_{\ell}} \cdot c_1(\omega_{\pi})^b\right)\in A^{\star}(B)\, .$$ When $\pi$ is the universal family of the moduli of polarized surfaces $$\pi: \mathcal S\to \mathcal M,$$ the $\kappa$-classes thus defined generate the tautological 
ring $\mathsf R^{\star}(\mathcal M)$ analogous to the well-studied tautological ring $\mathsf R^{\star}(\mathcal M_g)$ of the moduli of curves. 
Finding relations between the $\kappa$-classes in $\mathsf R^{\star}(\mathcal M)$ is a very interesting problem. 

In the case of the moduli of $K3$ surfaces, a strategy for $\kappa$-relations was laid out in \cite {mop} via the study of the virtual class of the Quot scheme; a different approach via Gromov-Witten theory was pursued in \cite{PY}. The discussion however makes sense for arbitrary polarized surfaces as well. In this approach, the center stage is taken by the calculation of the push-forwards $$\sum_{n=0}^{\infty} q^n \pi^{[n]}_{\star} \left(s_{i_1} (L_1^{[n]})\cdots s_{i_\ell}(L_{\ell}^{[n]})\right)\in A^{\star}(\mathcal M)$$ in terms of the classes $\kappa[a_1, \ldots, a_\ell, b].$ In the $K$-trivial case, or more ambitiously for arbitrary surfaces, it becomes important to obtain explicit formulas, thus generalizing the results of this paper.

\subsection {\bf Acknowledgements.} We are grateful for conversations and correspondence with N. Arbesfeld, D. Johnson, M. Kool, D. Maulik, A. Mellit, A. Okounkov, A. Szenes, and C. Voisin.

 A.M. was supported by the NSF through grant DMS 1601605. She thanks D. Maulik and the MIT mathematics department for their hospitality during the Fall of 2017. She also thanks MSRI and the organizers of the program {\it ``Enumerative Geometry Beyond Numbers"} in the Spring of 2018 for the wonderful conditions and inspiring atmosphere. D. O. was supported by the NSF through grants  DMS 1150675 and DMS 1802228. R.P. was supported by the Swiss National Science Foundation and
the European Research Council through
grants SNF-200020-162928, SNF-200020-182181,  ERC-2012-AdG-320368-MCSK,
ERC-2017-AdG-786580-MACI, SwissMAP, and the Einstein Stiftung. This project has received funding from the European Research
Council (ERC) under the European Union Horizon 2020 Research and
Innovation Program, grant agreement No. 786580.

\section {$K3$ surfaces} 
\subsection {Residue calculations} 
\label{rescal}
Our first goal is to prove Theorem \ref{t1}. We begin by explaining how the special formulas of Theorem \ref{p1} are predicted by the series in Theorem \ref{t1}. Conversely, we will prove in Section \ref{latelate} that these predictions are equivalent to the statement of Theorem \ref{t1}. 

We keep the same notation as in the introduction. Let $S$ be a $K3$ surface.
For a vector bundle $V\to S$ with Mukai vector $$v=\text{ch }(V)\sqrt{\text{td}(S)},$$ let 
 $$\chi = \chi (S,V)\, , \ \ c_1=c_1(V)\, , \ \   c_2 = c_2 (V)\, ,$$
 as in Section \ref{sttg}. Recall that $s=\text{rk }V$ and $r=s+1$. Taking the $K$-theory class $\alpha$ to be $V$,
the statement of Theorem \ref{t1} becomes 
$$\mathsf {S}(z)= 
\left (1+ (1+r)t \right )^{\left[(r-1)c_2+\left(-\frac{r}{2}+1\right)c_1^2-r^2+2r\right]} \cdot \left (1+tr \right )^{\left[-rc_2+\frac{r-1}{2}c_1^2+(r^2-1)\right]} \cdot \frac{1}{1 + r(1+r) t}\, .$$
For convenience, we define
\begin{equation}
 d = (r-1)c_2+\left(-\frac{r}{2}+1\right)c_1^2-r^2+2r\, .
 \label{a}
 \end{equation}
We then need to prove
 \begin{equation}\label{seg}\mathsf{S}(z) = \left (1+ (1+r)t\right )^{d} \cdot \left (1+t r\right )^{-d +\chi +1} \cdot \frac{1}{1 + r(1+r) t}.\end{equation}
An important observation is that the quantity \eqref{a}  is half 
 the dimension of the moduli space $\mathfrak M_v$ of stable sheaves of type $v$ on the K3 surface $S$,
 $$\dim \mathfrak M_v = \langle v, \, v \rangle + 2 = 2d\, .$$ 
We can express the individual Segre integrals as the residues $$\int_{S^{[n]}} s_{2n}  (V^{[n]}) = \frac{1}{2\pi i} \oint \mathsf{S} (z) \cdot \frac{dz}{z^{n+1}}\, .$$
Substituting equation \eqref{seg} in the above residue and using the change of variables \eqref{changez} $$z=t(1+rt)^r,$$ we are equivalently seeking to prove the
following formula for the top Segre classes:
\begin{eqnarray*}
\int_{S^{[n]}}s_{2n} (V^{[n]}) &=& \frac{1}{2\pi i} \oint \left ( 1 + (1+r)t \right )^{d} \cdot \left (1+rt \right )^{-d +\chi -rn} \cdot \frac{dt}{t^{n+1}}\\
& = & \text{Coeff}_{\,t^n} \left [ \left ( 1+ (1+r)t \right )^{d} \cdot \left (1+rt \right )^{-d +\chi -rn} \right]\, .
\end{eqnarray*}
The formula yields a remarkably simple answer in the following two cases, leading to the statement of Theorem \ref{p1}.
\begin{enumerate}
\item[(i)] 
When $d= 0,$ we expect
$$\int_{S^{[n]}}s_{2n}  (V^{[n]}) = r^n \binom{\chi - rn}{n}\, .$$
In particular, when $r = 2$ we recover the known line bundle result 
of \cite{mop},
$$\int_{S^{[n]}}s_{2n} (L^{[n]}) = 2^n \binom{\chi - 2n}{n}\, .$$

\item [(ii)] 
When $d = 1$, we expect
$$\int_{S^{[n]}}s_{2n}  (V^{[n]}) = r^n \left ( -r+\frac{1}{r} +\frac{\chi}{n} \right ) \binom{\chi-rn-1}{n-1}\, .$$
\end{enumerate}
\vskip.1in

 \subsection {\bf Vanishing results.}\label{vanres} We now prove the formulas above, thus establishing Theorem \ref{p1}. The proof is guided by the simple form of the expected answers. 
 
Throughout this Section, we let $S$ be a $K3$ surface of Picard rank one, $$\text{Pic} \, (S) =\mathbb Z H.$$ Let 
 $V\to S$ be an $H$-stable {\footnote {$H$-stability is Gieseker stability with respect to the
polarization $H$.}
vector bundle of rank $s>1$ with $c_1 (V) = H$, so that the Mukai vector equals $$v=(s, H, \chi-s).$$ We assume that $$\langle v, \, v \rangle = -2 \text{ or } \langle v, \, v \rangle = 0.$$ 
In both of
these cases, we know the moduli space of $H$-stable sheaves with Mukai vector $v$ is nonempty, see \cite{Y3}, Theorem $0.1$ for a general result. The moduli space is either a point, or a $K3$ surface, see Theorems $1.4$ and $3.6$ in \cite{Mu}. 

Since $s>1$, the locus of locally free sheaves $V$ in  
the moduli space with Mukai vector $v$ is {\it nonempty}.
Nonemptiness is
 obtained by invoking Yoshioka's classification of Mukai vectors yielding moduli consisting entirely of nonlocally free sheaves in Proposition $0.5$ of \cite {Y3}. His classification does not include vectors $v$ as above. 

Under these assumptions, we show the following vanishings.

\begin{proposition} We have
\begin{itemize} \label{p2}
\item [(i)] If $\langle v, \, v \rangle = -2,$ then $s_{2n} (V^{[n]}) = 0$ for $rn\,  \leq \, \chi (V)\,  < \, (r+1) n.$ 
\item [(ii)] If $\langle v, \, v \rangle = 0,$ then $s_{2n} (V^{[n]}) = 0$ for $rn +1 \leq \chi (V) < (r+1) n.$ 
\end {itemize}
\end{proposition}

\noindent {\it Proof.}  Recall that $V$ is said {\it $(n-1)$-very ample} if the evaluation map \begin{equation}\label{eval}H^0 (S, V) \to H^0 (S, V \otimes {\mathcal O}_Z)\end{equation} is surjective for all $Z \in S^{[n]}.$ This is equivalent to the surjectivity of the natural vector bundle map 
$$H^0 (S, V) \otimes \mathcal O_{S^{[n]}} \longrightarrow V^{[n]}\, \, \, \text{on} \, \, \, S^{[n]}.$$ We will show that this is the case under our numerics in Proposition \ref{p3} below. 

By the definition of Segre classes, whenever the evaluation map is surjective, we have the vanishing  $$s_j (V^{[n]} )= 0 \, \, \, \text{for} \, \, \,  j> h^0 (S, V) - \text{rank} \, V^{[n]}\, . $$
In our situation $$h^0(S, V)=\chi(V).$$ Indeed, by Serre duality and stability, $$
h^2(S, V)=h^0(S, V^{\vee})=0\, .$$ Furthermore, for $\chi>0$ we have 
\begin{equation}
h^1(S, V)=0\, . \label{qq78}
\end{equation}
The vanishing \eqref{qq78} is established
 in the proof of Proposition \ref{p3} below, where it is argued
that there do not exist nontrivial extensions 
$$0\to V^{\vee}\to E\to \mathcal O_S\to 0\, .$$ Therefore, $(n-1)$-very ampleness of $V$ implies the vanishing 
$$s_j  (V^{[n]}  )= 0 \, \, \, \text{for} \, \,  j > \chi (V) - sn\, .$$
In particular, we obtain $$s_{2n}  (V^{[n]}  )= 0\, \text { for } (r+1)n>\chi(V).$$\qed
\begin{proposition}\label{p3}
We have
\begin{itemize} 
\item [(i)] If $\langle v, \, v \rangle = -2$ and $\chi (V) \geq rn,$ then $V$ is $(n-1)$-very ample. 
\item [(ii)] If $\langle v, \, v \rangle = 0$ and $\chi (V) \geq r n + 1,$ then $V$ is $(n-1)$-very ample.
\end {itemize}
\end{proposition}

\noindent {\it Proof.} The evaluation map \eqref{eval} is surjective if $H^1 (V \otimes I_Z) = 0$ for all $Z \in S^{[n]}.$ Assume, for contradiction, that 
\begin{equation}\label{f699}
H^1(V\otimes I_Z)\neq 0\, .
\end{equation}
The cohomology group \eqref{f699}
is Serre dual to $\text{Ext}^1(I_Z, V^{\vee})$ which is the space of extensions
\begin{equation}\label{ex}0 \to V^{\vee} \to E \to I_Z \to 0\, .\end{equation}
By assumption, we have a non-split extension. 

If $V$ is $H$-stable with $c_1 (V) = H,$ then Lemma 2.1 of \cite{Y2} 
shows that the middle term $E$ of any non-split extension
$$0 \to V^{\vee} \to E \to I_Z \to 0$$ is $H$-stable. 

For the benefit of the reader, let us recall the argument in \cite {Y2} in our context. Since  $c_1(V) = H$ is primitive, note first that the $H$-stability of $V$ implies that $V$ is in fact slope-stable, therefore $V^{\vee}$ is slope-stable, hence $V^{\vee}$ is $H$-stable as well. 

Assuming now that $E$ is not $H$-stable, let $G \hookrightarrow E$ be the maximal semistable destabilizing subsheaf. We have $${\text{rk}\, G} < {\text{rk}\, E} \, \, \, \text{and} \, \,  \,   \mu_G \geq \mu_E > 
\mu_{V^{\vee}}\, .$$ We see that $G$ cannot be a subsheaf of the kernel $V^{\vee}$ since this would contradict the $H$-stability of the latter. Therefore, we have a nonzero morphism $$\phi: G \to I_Z\, .$$
The $H$-semistability of $G$ now gives  $\mu_G \leq 0\, .$ Writing $c_1(G) = a H$ we deduce $a\leq 0$. Since $\mu_G \geq \mu_E$, we obtain $$\frac{a H^2}{\text{rk}\, G} \geq -\frac{H^2}{\text{rk}\, E} \implies a>-1\, .$$ Therefore, we have $$a = 0 \text{ and }c_1(G) = 0\, .$$ 
Furthermore, if nonzero, the kernel $K$ of $\phi$ is a subsheaf of $V^{\vee}$ of slope greater than or equal to zero, contradicting the $H$-stability of $V^{\vee}.$ Thus, we have  $$\phi: G \to I_Z \text { is injective} \implies G = I_W$$ for a zero-dimensional subscheme $W \subset S$.
We consider the exact sequence $$0 \to G \to I_Z \to Q \to 0$$ and the associated sequence of extension groups
$$\text{Ext}^1 (Q, V^{\vee}) \to \text{Ext}^1 (I_Z, V^{\vee}) \stackrel{f}{\rightarrow}\text{Ext}^1 (G, V^{\vee})\, .$$ The first group is zero by Serre duality
since $Q$ is supported at finitely many points. We conclude that $f$ is 
injective. However, the image of the extension \eqref{ex} in $\text{Ext}^1 (G, V^{\vee})$ is trivial. The contradiction shows that $E$ must be $H$-stable.

Now, we calculate
\begin{eqnarray*}
\chi (E, E) &= &\chi (V^{\vee}, V^{\vee}) + 2 \, \chi \left (V^{\vee}, \, I_Z \right ) + \chi (I_Z, \, I_Z)\\ & =& - \langle v, \, v \rangle + 2 (\chi (V) - sn) + 2 - 2n \\ &= &2 \left ( - \frac{\langle v, \, v \rangle}{2} + \chi (V) - rn + 1 \right).
\end{eqnarray*}
In both cases of Proposition \ref{p3}, we obtain 
$$\chi (E, E)  \geq 4\implies \text{ext}^0(E, E)+\text{ext}^2(E, E)\geq 4\, .$$ Therefore, by Serre duality, $$\text{ext}^0(E, E) \geq 2\, .$$ 
Since $E$ is $H$-stable and therefore simple, we have a contradiction. Thus,  \eqref{f699}
does not hold.
The proofs of Proposition \ref{p3} and Proposition \ref{p2}
are therefore complete. \qed 
\vskip.1in

\subsection{Proof of Theorem \ref{p1}} We consider case (i) of the Theorem.  Fix the rank $s>1$ throughout. \footnote{
For $s=1$, the functions $A_1, A_2$ are known by \cite{mop,voisin}. All three functions $A_0, A_1, A_2$ for $s=1$ follow from the case $s>1$ via the polynomiality argument of Section \ref{cc}.} By Theorem $4.1$ in \cite{EGL}, over $K3$ surfaces, the Segre integral $$\int_{S^{[n]}} s_{2n}(V^{[n]})$$ is given by a universal polynomial of degree $n$ in $c_1(V)^2$ and $c_2(V)$. When $$\langle v, v\rangle =-2$$ both $c_1(V)^2$ and $c_2(V)$ can be expressed in terms of $\chi=\chi(V)$, and therefore 
$$\int_{S^{[n]}} s_{2n}(V^{[n]})=\mathsf P_n(\chi),$$ for a degree $n$ universal polynomial $\mathsf P_n$.

For each $\chi>s$, let $(S, H)$ be a $K3$ surface of Picard rank one $\text{Pic} (S)=\mathbb Z H$ and of genus $$H^2=2g-2,\quad g=s(\chi-s).$$ Let $V$ be the unique $H$-stable vector bundle with Mukai vector $v=(s, H, \chi-s)$ whose existence we noted in the beginning of Section \ref{vanres}. Since $\langle v, v\rangle =-2$, 
Proposition \ref{p2} (i) shows that the top Segre class of $V^{[n]}$ vanishes when $$rn\leq \chi<(r+1)n.$$ These $n$ values of $\chi$ are the $n$ roots of the polynomial $\mathsf P_n$. 
Therefore, $$\mathsf P_n(\chi)=c \cdot \binom{\chi-rn}{n}\implies \int_{S^{[n]}} s_{2n}(V^{[n]})=c \cdot \binom{\chi-rn}{n}$$ for some constant $c$. We will identify the constant
$c=r^n$ by proving 
 \[\int_{S^{[n]}} s_{2n}(V^{[n]})=\frac{r^n}{n!}\chi^n+\text{lower order terms (l.o.t.) in }\chi.\]

The argument is most naturally expressed by
rewriting \eqref{eq1} in exponential form,
 $$\sum_{n=0}^{\infty} z^n \int_{S^{[n]}} s_{2n}(V^{[n]})=
\exp \Big(\bar {A}_0(z) \cdot c_1(V)^2+\bar {A}_1(z) \cdot c_2(V) +\bar{A}_2(z)\Big)\,.$$ 
Restricting to the line $\langle v, v\rangle =-2$ in the $(c_1^2, c_2)$-plane,
we obtain 
$$\sum_{n=0}^{\infty} z^n \int_{S^{[n]}} s_{2n}(V^{[n]})=
\exp \Big(U(z) \cdot \chi +T(z)\Big)\, ,$$ for power series $U$ and $T$. Let $$U(z)=u_1z+u_2z^2+\ldots .$$ Extracting the coefficient of $z^n$ in the above expression yields $$\int_{S^{[n]}} s_{2n}(V^{[n]})=\frac{u_1^n}{n!} \chi^n+\text{ l.o.t.}\ .$$ In particular, for $n=1$, we obtain $$\int_S s_2(V)=u_1 \chi+\text{ l.o.t.}\ .$$ Direct calculation shows $s_2(V)=r\chi-r^2$ so that $u_1=r.$ 
After substitution,
we find the leading term to be $$\frac{u_1^n}{n!}=\frac{r^n}{n!}\, . $$
  
For part (ii) of the Theorem, Proposition \ref{p2} (ii) gives only $n-1$ roots of the polynomial expressing the Segre integral. Combined with the leading term calculations, 
we conclude that if $\langle v, v\rangle =0$ we have $$\int_{S^{[n]}} s_{2n}(V^{[n]})=\frac{r^n}{n}(\chi+c)\binom{\chi-rn-1}{n-1}\, ,$$ for some constant $c$. We will prove
 $$c=n\left(-r+\frac{1}{r}\right)$$ 
by computing the next term in the Segre polynomial in $\chi$. Indeed, for the right hand side, the $\chi^{n-1}$ coefficient is easily seen to be \begin{equation}\label{coeff1}\frac{r^n}{n} \left(\frac{c}{(n-1)!}-\frac{n\left(r+\frac{1}{2}\right)}{(n-2)!}\right).\end{equation}

We compute the $\chi^{n-1}$ coefficient on the left. 
As before, we have $$\sum_{n=0}^{\infty} z^n \int_{S^{[n]}} s_{2n}(V^{[n]})=
\exp \Big(U(z) \cdot \chi +T(z)\Big)$$ 
where $$U(z)=u_1z+u_2z^2+\ldots\,, \ \ \  T(z)=t_1z+t_2z^2+\ldots\, .$$
We obtain
$$\int_{S^{[n]}} s_{2n}(V^{[n]})=\frac{u_1^n}{n!}\cdot \chi^{n}+\left(\frac{u_1^{n-1}t_1}{(n-1)!}+\frac{u_1^{n-2}u_2}{(n-2)!}\right)\cdot \chi^{n-1}+\ldots\, .$$ 
When $n=1$, we find
 $$\chi r- r^2+1=\int_S s_2(V)=u_1 \chi + v_1\implies u_1=r, \,\, t_1=-r^2+1\,.$$ When $n=2$, for isotropic vectors $\langle v, v\rangle=0$, 
equation \eqref{2pt} simplifies to $$\int_{S^{[2]}} s_4(V^{[2]})=r^2\left(-r+\frac{1}{r}+\frac{\chi}{2}\right)\cdot (\chi-2r-1)\, .$$ 
Using the above asymptotics, we have $$ \int_{S^{[2]}} s_{4}(V^{[2]})=\frac{u_1^2}{2}\cdot \chi^{2}+\left(u_1t_1+u_2\right)\cdot \chi+\ldots\implies u_2=-r^3-\frac{r^2}{2}.$$
Thus, for arbitrary $n$, the $\chi^{n-1}$-coefficient equals \begin{equation}\label{secondeq}\frac{r^{n-1}(1-r^2)}{(n-1)!}+\frac{r^{n-2}(-r^3-\frac{r^2}{2})}{(n-2)!}.\end{equation} 

Comparison of \eqref{coeff1} and \eqref{secondeq} yields the requisite value for the constant $c$, completing the argument. \qed

\subsection {Proof of Theorem \ref{t1}} \label{latelate}
By \cite{EGL}, the generating series of Segre integrals takes the form 
$$\mathsf{S} (z) = \sum_{n=0}^{\infty} z^n \int_{S^{[n]}} s_{2n}(V^{[n]})=a_0(z)^{c_2 (V)} \cdot a_1 (z)^{c_1(V)^2} \cdot a_2 (z)^{\chi(\mathcal O_S)}\, ,$$ 
for three power series $a_0 (z), a_1 (z), a_2 (z) \in \mathbb Q [[z]].$ 
We show that these power series are as claimed by Theorem \ref{t1}: 
$$a_0(z)=A_0(z), \,\,\,a_1(z)=A_1(z), \,\,\,a_2(z)=A_2(z).$$

Fix an integer $\chi>s$, and let $g=s(\chi-s).$ Consider $(X, H)$ a $K3$ surface of genus $g$,
 and let $V$ be a vector bundle with Mukai vector $$v=(s, H, \chi-s)\implies \langle v, v\rangle=-2\, .$$ 
We calculate $$c_2(V)=\chi(s-1)-s^2+2s-1,\,\, 
c_1(V)^2=2(s\chi-1-s^2)\, .$$ The residue calculations of Section \ref{rescal} and the first part of Theorem \ref{p1} together  
imply $$\left(A_0^{s-1} A_1^{2s}\right)^{\chi} \left(A_0^{-s^2+2s-1} A_1^{-2-2s^2}A_2^2\right)=\left(a_0^{s-1} a_1^{2s}\right)^{\chi} \left(a_0^{-s^2+2s-1} a_1^{-2-2s^2}a_2^2\right).$$ 
Since $\chi$ is arbitrary, we obtain \begin{equation}\label{e1}A_0^{s-1} A_1^{2s}=a_0^{s-1} a_1^{2s}\, ,\end{equation} \begin{equation}\label{e2}A_0^{-s^2+2s-1} A_1^{-2-2s^2}A_2^2=a_0^{-s^2+2s-1} a_1^{-2-2s^2}a_2^2\, .\end{equation}

We now derive an additional equation using isotropic bundles corresponding to the second part of Theorem \ref{p1}.  To this end, let $g=s(\chi-s)+1$. Let 
$V\to S$ be a vector bundle with Mukai vector 
$$v=(s, H, \chi-s)\implies \langle v, v \rangle=0\, .$$
 We then calculate
 $$c_2(V)=\chi(s-1)-s^2+2s\, , \ \ \  c_1(V)^2=2(s\chi-s^2)\, .$$ 
Repeating the above argument for the new numerics, we replicate equation \eqref{e1} and in addition we obtain
\begin{equation} \label {e3}A_0^{-s^2+2s} A_1^{-2s^2}A_2^2=a_0^{-s^2+2s} a_1^{-2s^2}a_2^2\, .\end{equation} Constraints \eqref{e1}, \eqref{e2}, \eqref{e3} show that $A_i=a_i \text { for } 1\leq i\leq 3,$ as claimed.

\subsection {Chern classes.} \label{cc} 
Let $\alpha\to S$ be a $K$-theory class of rank $s$ on a $K3$ surface $S$, and let $s=r+1$. 
As a corollary of Theorem \ref{t1}, we also obtain
expressions for the Chern classes of tautological bundles,
 $$\sum_{n=0} z^n \int_{S^{[n]}} c_{2n}(\alpha^{[n]})=\tilde A_0(z)^{c_2(\alpha)}\cdot \tilde A_1(z)^{c_1^2(\alpha)}\cdot \tilde A_2(z)^{\chi(\mathcal O_S)}\, .$$ 
After the change of variables $$z=t(1-rt)^{-r}\, ,$$
the series are given by 
 \begin{eqnarray*}\tilde A_0(z)&=&(1-rt)^{-r}\cdot (1+(-r+1)t)^{r+1}\, , \\\tilde A_1(z)&=& (1-rt)^{\frac{r-1}{2}} \cdot (1+(-r+1)t)^{-\frac{r}{2}}\, , \\\tilde A_2(z)&=&(1-rt+r^2t)^{-\frac{1}{2}}\cdot (1-rt)^{\frac{r^2-1}{2}} \cdot (1+t(-r+1))^{-\frac{r^2}{2}-r}.\end{eqnarray*}

\begin{remark} The rank $s=2$ case (corresponding to $r=1$) is
easily calculated from geomety by assuming the
existence of a transverse section of $V$, see \cite{J}. 
The derivation using our formulas is also simple,
$$\tilde A_0 = \frac{1}{1-t}\,, \ \ \ \tilde A_1= 1\, ,\ \ \ \tilde A_2=1\,, \ \ \ z=  \frac{t}{1-t}\implies \tilde A_0=1+z.$$
We find 
\begin{equation}\label{rank2222}\sum_{n=0} z^n \int_{S^{[n]}} c_{2n}(V^{[n]}) = (1+z)^{c_2(V)}\, ,\end{equation}
which is the correct answer.
\end{remark}

The evaluation of $\tilde A_0$, $\tilde A_1$, and $\tilde A_2$ for
Chern classes follows by
 regarding the Segre integrals of Theorem \ref{t1} as functions on the $K$-theory of the surface $S$ which depend polynomially on
 $$\text{rank }\alpha=s,\, \ \ c_1(\alpha)^2\, , \ \ c_2(\alpha)\, ,$$ 
see \cite {EGL}. 
Having established these polynomials for positive $s>0$, we may then also allow $s$ to be negative, and replace $\alpha$ by $-\alpha$.

As before, if $V$ is a bundle with $\langle v, v\rangle=-2$, we obtain 
$$\int_{S^{[n]}}c_{2n}(V^{[n]})=(-r)^{n} \binom{-\chi+rn}{n}\, .$$ 
In particular, paralleling Proposition \ref{p2}, when $\langle v, v\rangle=-2$ we obtain the vanishing $$c_{2n} (V^{[n]}) = 0  \text { for }(s-2)n<\, \chi (V)\, \leq (s-1) n.\, $$

\vskip.1in

\section{Abelian and Enriques surfaces}

\subsection {\bf Abelian surfaces.} Having established the $K3$ case, we can also determine the Segre integrals over abelian or bielliptic surfaces $S$. Indeed, by the universality results of \cite {EGL}, we have 
$$\sum_{n=0}^{\infty} z^n \int_{S^{[n]}}  s_{2n} (V^{[n]}) = A_0(z)^{c_2(V)}\cdot A_1 (z)^{c_1^2 (V)}.$$ Residue calculations give the coefficients of the right hand side $$ \int_{S^{[n]}}  s_{2n} (V^{[n]})= \text{Coeff}_{\,t^n} \left [ \left ( 1+ (1+r)t \right )^{d} \cdot \left (1+rt \right )^{-d +\chi -rn-1} \cdot \left(1+r(r+1)t\right)\right],$$ where $\langle v, v\rangle =2d$. In particular, we obtain
 the following generalization of the $s=1$ result of \cite {mop}.
 \begin{proposition} 
Let $V\to S$ be a vector bundle of rank $s=r-1$ on an abelian or
bielliptic surface $S$. 
If $\langle v, v\rangle=0$, then $$ \int_{S^{[n]}}  s_{2n} (V^{[n]})=r^n\cdot \frac{\chi}{n} \cdot \binom{\chi-rn-1}{n-1}\, .$$ \end{proposition} \noindent 

\subsection {\bf Enriques surfaces and strange duality.} \label{esdd} When $S$ is an Enriques surface, the link between the Segre and Verlinde series yields individual equalities of intersection numbers, as the numerical data can be suitably matched. The Proposition below parallels Conjecture $2.2$ of \cite {J} formulated for del Pezzo surfaces. 
\comment{The result is implicit in Theorem $1$, its precise matching with Conjecture $3.3$ of \cite{J}, and the calculations contained in the proof of Theorem 3.8 in \cite {J}. For the benefit of the reader, we also give a direct argument.} The corresponding result does not hold for other $K$-trivial surfaces. 
  
\begin{proposition} \label{ccqq} Let $V\to S$ be a vector bundle of rank $s=r+1$, determinant $L$, and $\chi(V)=(r-1)n+1.$ 
We have $$\int_{S^{[n]}} c_{2n}(V^{[n]})=\chi(S^{[n]}, L_{(n)}\otimes E^r)\, .$$
\end{proposition}

\comment{This result and the numerics under which it holds are both motivated by strange duality, as explained in 
\cite {J}. }The left hand side interprets enumeratively the Verlinde number on the right as an integral over the Hilbert scheme -- counting the expected (finite) number of quotients $V^{\vee}\to I_Z$ with $Z$ a subscheme of length $n$.  Implications of Proposition \ref{ccqq} for strange duality over $K3$ and Enriques surfaces will be taken up elsewhere. 
 
\proof Write $\chi=\chi(L)$ and note that the assumption of the Proposition translates into $$c_2(V)=\chi-(r-1)(n-1)\, .$$ 
By Subsection \ref{cc}, the Chern integral on the left hand side is the coefficient of $z^n$ in the series $$\tilde A_0(z)^{c_2(V)} \cdot \tilde A_1(z)^{c_1^2(V)} \cdot \tilde A_2(z)\, ,$$ or equivalently,
 the residue of the differential form 
$$\tilde A_0(z)^{\chi-(r-1)(n-1)}\cdot \tilde A_1(z)^{2\chi-2} \cdot \tilde A_2(z) \cdot \frac{dz}{z^{n+1}}\, .$$ 
Using the change of variables $z=t(1-rt)^{-r}$, the differential form becomes \begin{equation}\label{form1} (1-tr(1-r))^{\frac{1}{2}}\cdot (1-tr)^{-\chi+r^2\left(n-\frac{1}{2}\right)-\frac{1}{2}}\cdot \left(1+\left(1-r\right)t\right)^{\chi-r^2\left(n-\frac{1}{2}\right)+(n-1)}\cdot \frac{dt}{t^{n+1}}\end{equation} For the right hand side, we use the calculations of Lemma 5.2 of \cite {EGL}: $$\sum w^n\chi(S^{[n]}, L_{(n)}\otimes E^r)=F(w)^{\frac{1}{2}}\cdot G(w)^{\chi}$$ where $$F(w)=\frac{(1+u)^{r^2}}{1+r^2u}\,, \ \ \ G(w)=1+u$$ for the change of variables $$w=u(1+u)^{r^2-1}.$$
Therefore, $\chi(S^{[n]}, L_{(n)}\otimes E^r)$ is the residue of the differential form \begin{equation}\label{form2}F(w)^{\frac{1}{2}}\cdot G(w)^{\chi} \cdot \frac{dw}{w^{n+1}} =(1+r^2u)^{\frac{1}{2}}\cdot (1+u)^{\chi-r^2\left(n-\frac{1}{2}\right)+(n-1)}\cdot \frac{du}{u^{n+1}}\end{equation} The change of variables $$u=\frac{t}{1-tr}$$ matches the two differential forms \eqref{form1} and \eqref{form2} and completes the proof.\qed

\section {$K3$ blowups}  

In order to determine the functions $A_3, A_4$, we need to consider surfaces which are not $K$-trivial. We will look at two different families of examples over the blowup of a $K3$ at a point. These examples are stated in Theorem \ref{t4}. We will prove this theorem first. As before, Reider-type arguments play a key role in the calculation. In addition, several excess intersection calculations are needed. Afterwards we show Theorem \ref{conj2}: the integrals calculated in Theorem \ref{t4} give the rank $2$ series for all surfaces. 

\subsection {Stability in extensions} Let $$\pi: S \to X$$ be the blowup of a $K3$ surface $X$ at one point. We assume that $\text{Pic} \, X = \mathbb Z H.$ Then $$\text{Pic}\, S = \mathbb Z H + \mathbb Z E,$$ where $E$ is the exceptional divisor on $S$ and $H$ denotes the pullback to $S$ of the ample Picard generator on $X$. Note that $H$ is a nef line bundle on $S$. 

\vskip.1in
The notion of $H$-stability and $H$-semistability have the usual meaning:
\begin{definition} We say that a torsion-free sheaf $V$ on $S$ is $H$-stable if for any nonzero subsheaf $G$ of strictly smaller rank, we have  $$\frac{c_1 (G) \cdot H}{\text{rk} \, G} < \frac{c_1 (V) \cdot H}{\text{rk}\, V}.$$
For the notion of $H$-semistability the inequality is not required to be strict. 
\end{definition}

\vskip.1in

We first prove several statements regarding vanishing of cohomology which will be useful in the proof of Theorem \ref{t4}. We consider a vector bundle $V_0\to X$ with Mukai vector $v_0$ such that \begin{equation}\label{apt}\text{rank }V_0=2,\,\,\, c_1 (V_0) = H, \,\,\, \langle v_0, v_0\rangle =-2, \,\,\,V_0\,\, \,\text{is Gieseker }H\text{-stable}.\end{equation} The existence of the bundle $V_0$ was noted in the beginning of Section \ref{vanres}: $V_0$ is the unique point in the moduli space of Gieseker $H$-stable bundles with Mukai vector $v_0$. Such vector bundles are necessarily {\it rigid}, that is $$\text{Ext}^1(V_0, V_0)=0.$$
Let \begin{equation}\label{bunv}V = \pi^{\star} V_0 \otimes E^{-k} \, \, \, \text{on} \, \, \,  S, \, \, \, \text{for} \, \, \, k \geq 0.\end{equation}
In order to compute Segre classes of $V^{[n]}\to S^{[n]}$, we  analyze the surjectivity of the evaluation map:
\begin{equation}
\label{ev}
H^0 (V) \otimes \mathcal O_{S^{[n]}} \to V^{[n]} \, \, \text{on} \, \, S^{[n]}.
\end{equation} 
Thus, for a zero dimensional subscheme $Z$ of length $n$, we investigate the vanishing of $$H^1(S, V\otimes I_Z)=\text{Ext}^1(V^{\vee}, I_Z)=\text{Ext}^1(I_Z, V^{\vee}\otimes E)^{\vee}.$$ Assume that a non-zero element exists in the extension group. We prove 
\begin{lemma}\label{stabb} \begin{itemize} 
\item [(i)] Assume $k\geq n-1$. 
For any nonsplit extension \begin{equation}\label{defext}0 \to V^{\vee} \otimes E \to W \to I_Z \to 0\end{equation} with $Z \in S^{[n]},$ the middle sheaf $W$ is $H$-stable. 
\item [(ii)] Assume $k=n-2$. The same conclusion holds, unless the scheme $Z$ is supported on the exceptional divisor $E$. \end{itemize} 
\end{lemma}

\proof  Assume $G \hookrightarrow W$ is an $H$-semistable destabilizing subsheaf. We then have $$\mu_G \geq \mu_W > \mu_{V^{\vee} \otimes E}.$$ 

Since $V_0$ is $H$-slope stable on $X$, $V_0^{\vee}$ is also $H$-slope stable. The pullback $\pi^{\star} V_0^{\vee}$ is $H$-stable on $S$. Thus $V^{\vee}\otimes E$ is also $H$-stable. As a consequence, we see that $G$ cannot be a subsheaf of the kernel $V^{\vee} \otimes E$ since this would contradict the $H$-stability of the latter. 

Thus we have a nonzero morphism $\phi: G \to I_Z.$
The $H$-semistability of $G$ now gives  $\mu_G \leq 0.$ If $c_1(G) = a H + b E,$ we deduce
$$0\geq \mu_G\geq \mu_{W}\implies 0 \geq \frac{a H^2}{\text{rk}\, G} \geq -\frac{H^2}{\text{rk}\, W},$$ from which we conclude that $a = 0$ and $\mu_G = 0.$ Furthermore, if nonzero, the kernel $K$ of $\phi$ is a subsheaf of $V^{\vee} \otimes I_Z$ of slope greater than or equal to zero, contradicting the $H$-semistability of $V^{\vee} \otimes E.$ Thus we have in fact that $\phi: G \to I_Z$ is injective, so $$G = I_U \otimes E^{-m}$$ for a zero-dimensional subscheme $U \subset S$ and $m \geq 0.$ In particular $Z$ is supported on the
exceptional curve $E$ and on $U$.
We turn to the exact sequence $$0 \to G \to I_Z \to Q \to 0$$ and the associated sequence of extension groups
$$\text{Ext}^1 (Q, V^{\vee} \otimes E) \to \text{Ext}^1 (I_Z, V^{\vee} \otimes E) \stackrel{\alpha}{\longrightarrow} \text{Ext}^1 (G, V^{\vee} \otimes E).$$ The image of the extension \eqref{defext} $$0\neq e\in \text{Ext}^1 (I_Z, V^{\vee} \otimes E)$$ under $\alpha$ is trivial since the resulting extension is seen to be split. Turning to the first extension group, we have
$$\text{Ext}^1(Q, V^{\vee}\otimes E)^{\vee}= H^1 (V \otimes Q).$$ Recall the notation $V=\pi^{\star} V_0\otimes E^{-k}$. We claim that $$H^1(V\otimes Q)=0, \text{ for } k\geq n-1.$$ This would imply $\alpha$ is injective, and in turn that the original extension $e$ splits -- a contradiction. 

Note that $Q$ is supported on the exceptional divisor and a finite number of points
in $S.$ To prove the claimed vanishing, let us first assume that $U$ is empty. The argument is best illustrated by the case $m=1$. In this case, the defining sequence $$0\to E^{-1}\to I_Z\to Q\to 0$$ shows that $Q=I_{Z/E}$ is the ideal sheaf of $n$ points on the exceptional divisor. Since $$V|_{E}= \mathbb C^2\otimes \mathcal O_E(k),$$ we see that $$H^1(V\otimes Q)= H^1(\mathcal O_E(k-n))\otimes \mathbb C^2=0 \text{ for } k\geq n-1.$$ We now consider the case of arbitrary $m$, continuing to assume that $U=\emptyset$. We let $W_\ell=\ell E$ be the scheme defined by the ideal sheaf $E^{-\ell}$. Set $Z_m=Z$ and $Q_m=Q$. The defining exact sequence $$0\to E^{-m} \to I_Z\to Q\to 0$$ shows that $$0\to Q_{m}\to \mathcal O_{W_{m}}\to \mathcal O_{Z_m}\to 0.$$ Inductively define the scheme-theoretic intersection $$Z_\ell=Z_{\ell+1}\cap W_\ell \hookrightarrow W_{\ell+1}$$ for $1\leq \ell \leq m-1,$ and let $Q_\ell$ be given by the exact sequence $$0\to Q_{\ell}\to \mathcal O_{W_{\ell}}\to \mathcal O_{Z_{\ell}}\to 0.$$
Clearly, $$\text{length }(Z_{\ell})\leq \text {length } (Z_m)=n.$$ We will show inductively that \begin{equation}\label{va}H^1(V\otimes Q_{\ell})=0\end{equation} for all $\ell\leq m$. The base case $\ell=1$ was verified above. For the inductive step, form the diagram 
\begin{center}
$\xymatrix{
& \ar[d] 0 & \ar[d]0 & \ar[d]0 &\\
0\ar[r] & K_\ell \ar[r]\ar[d] & \mathcal O_E(\ell-1)\ar[r]\ar[d] & M_\ell\ar[r]\ar[d] & 0 \\
0\ar[r] & Q_\ell \ar[r]\ar[d] & \mathcal O_{W_\ell}\ar[r]\ar[d] &\mathcal O_{Z_\ell}\ar[r]\ar[d] & 0 \\
0\ar[r] & Q_{\ell-1} \ar[r]\ar[d] & \mathcal O_{W_{\ell-1}}\ar[r]\ar[d] & \mathcal O_{Z_{\ell-1}}\ar[r]\ar[d] & 0\\
& 0 & 0 & 0 & }$
\end{center} 
Note that the support of $M_{\ell}$ has length at most $n$. As we already noted $$V|_{E} = \mathbb C^2\otimes \mathcal O_E(k).$$ If $\ell\geq 1$ and $k+1\geq n$, we have just enough positivity to ensure that the morphism \begin{equation}\label{susur}H^0(V\otimes \mathcal O_E({\ell}-1))=\mathbb C^2\otimes H^0(E, \mathcal O_E(k+(\ell-1))\to H^0(V\otimes M_{\ell})\end{equation} is surjective. Using the first row we conclude $$H^1(V\otimes K_{\ell})=0.$$ Using the first column and invoking the inductive hypothesis, we conclude $$H^1(V\otimes Q_{\ell})=0.$$ This completes the argument when $U$ is empty. For the general case, let $T$ denote the scheme with ideal $E^{-m}\otimes I_U$, and recall $W_m$ had ideal $E^{-m}.$ The defining exact sequence gives $$0\to I_T\to I_Z\to Q\to 0 \implies 0\to Q\to \mathcal O_T\to \mathcal O_Z\to 0.$$ By composing the first map with the canonical restriction $\mathcal O_T\to \mathcal O_{W_m}$ we obtain an exact sequence \begin{equation}\label{new}0\to Q\to \mathcal O_{W_m} \to A\to 0.\end{equation} Furthermore $$\mathcal O_Z\to A\to 0,$$ so $A=\mathcal O_{\widetilde{Z}}$ is supported on at most $n$ points. Then the previous argument applied to the exact sequence \eqref{new}
$$0\to Q\to \mathcal O_{W_m}\to \mathcal O_{\widetilde Z}\to 0$$ 
gives the vanishing $H^1(V\otimes Q)=0$.  This finishes the proof when $k=n-1$. 

When $k=n-2$, the same argument carries through, unless $m=1$ and $Z$ is contained in $E$, as one can easily check going through the details, in particular by examining \eqref{susur}. \qed
\vskip.1in

\begin{lemma}\label{endlemma}
 If $W$ is $H$-stable, then for any $k \geq 0$,  the dimension of $\text {Hom}(W, W\otimes E^k)$ equals $1.$
\end{lemma}

\proof The $H$-stability of $W$ implies that any nonzero homomorphism $\phi: W \to W \otimes E^k$ is injective. This can be seen as usual by examining the kernel and image of $\phi$. Let $\lambda$ be an eigenvalue of $\phi$ at a point $p$ in $S$ where $W$ is locally free and $p\not \in E$. Let $I$ be the canonical homomorphism $W \to W \otimes E^k$ obtained by $W$-twisting the unique section $\mathcal O \to E^k.$ We claim $$\phi=\lambda I.$$ Indeed, assuming otherwise, set $$\psi:=\phi-\lambda I\neq 0.$$ By the first line of the proof, the morphism $\psi:W\to W\otimes E^{k}$ must be injective. 
Consider the induced morphism $$\det \psi: \det W\to \det W\otimes E^{kr}.$$ Writing as before $I$ for the $\det W$-twisting of the unique section of $E^{kr}$, we conclude $$\det \psi = \mu I$$ for some constant $\mu$. 
By construction, $\det \psi$ vanishes at $p$. However $I$ only vanishes along $E$ and $p\not \in E$. This shows that $\mu=0$, so $\det \psi = 0$. This contradicts the fact that $\psi$ is injective. 
\qed

\vskip.2in

\begin{lemma}\label{propexpl}
Let $V=\pi^{\star} V_0\otimes E^{-k}$ with $V_0$ a rank $2$ bundle satisfying \eqref{apt}. Assume \footnote{This assumption uniquely specifies the genus $H^2=2g-2$ of the $K3$ surface $X$ in terms of $n$ and $k$, as well as the Mukai vector $v_0$ for each such $K3$ surface.} $\chi(V)=4n-1$.
\begin{itemize} 
\item [(i)] If $k=n-2$ then $$H^1(V\otimes I_Z)=0$$ for all $Z\in S^{[n]}$ unless $Z\subset E$. 
\item [(ii)] If $k=n-1$ then any nontrivial extension $W$ in \eqref{defext} is a rigid sheaf. 
\item [(iii)] In both cases, $H^1(V)=H^2(V)=0.$ 
\end{itemize}
\end{lemma}

\proof We continue working with the exact sequence \eqref{defext}. If $Z\not \subset E$, $W$ must be $H$-stable by Lemma \ref{stabb}. Using Lemma \ref{endlemma}, we compute $$\chi(W, W)\leq \text{ext}^0(W, W)+\text {ext}^2(W, W)=\text{ext}^0(W, W)+\text{ext}^0(W, W\otimes E)=2.$$
On the other hand, from the defining exact sequence \eqref{defext}, we calculate:
\begin{eqnarray*}
\chi (W, W) &= &\chi (V^{\vee}, V^{\vee}) +  \chi \left (V^{\vee} \otimes E, \, I_Z \right ) +  \chi \left (V^{\vee},  \, I_Z \right )+  \chi (I_Z, \, I_Z)\\ & =& - \langle v_0, \, v_0 \rangle + \left ( \chi (V_0) - 2n - {(k+1)(k+2)}\right ) +  \left ( \chi (V_0) - 2n - k(k+1) \right )+  2 - 2n \\ &= &2 \left ( - \frac{\langle v_0, \, v_0 \rangle}{2} + \chi (V_0) - 3n + 1 - (k+1)^2\right)\\  & = & 2 \left ( - \frac{\langle v_0, \, v_0 \rangle}{2} + \chi (V) - 3n-k\right)=2\left(n-k\right).
\end{eqnarray*}
\noindent
Note now that if $k=n-2$, then $$\chi(W, W)=2(n-k)=4$$ which is a contradiction. This establishes (i).  

If $k=n-1$, then $\chi(W, W)=2$. Since $\text{Ext}^0(W, W)=\text{Ext}^2(W, W)=\mathbb C$ we find that $\text{Ext}^1(W, W)=0$ so $W$ is rigid. This establishes (ii). 

Note that the same argument for $Z=\emptyset$ shows that $H^1(V)=0$, while $H^2(V)=0$ follows by stability. \qed
\vskip.1in
We will analyze the two situations (i) and (ii) in Propositions \ref{asp}, \ref{sub}, and \ref{p7} below. 

\subsection {Excess calculations for $k=n-2$} The goal of this section is to prove the following result which corresponds to Theorem \ref{t4}\,(i). 
\begin{proposition} \label{asp} Let $S$ be the blowup of a $K3$ surface at one point. Let $$V=\pi^{\star} V_0\otimes E^{-(n-2)}$$ with $V_0$ a rank $2$ bundle satisfying \eqref{apt}. Assume furthermore that $\chi(V)=4n-1$. Then $$\int_{S^{[n]}} s(V^{[n]})=(-1)^n(2n+1).$$
\end {proposition} 

\proof We compute $s_{2n} (V^{[n]})$ as an excess intersection over the nonsurjectivity locus of the evaluation map \eqref{ev}. 
By Lemma \ref{propexpl}, this nonsurjectivity locus consists of those $Z$ with $Z\subset E$, or in other words $$Z\in E^{[n]}\simeq\mathbb P^n.$$
For the excess calculation, we need to interpret the Segre class $s_{2n}$ in connection with the top Chern class of a vector bundle, where the excess formula is easier to understand. This connection is as follows, cf. Fulton \cite {F}.  

\vskip.2in

Let $G (2n, 4n-1)$ be the Grassmannian of $2n$ planes $\Lambda\hookrightarrow H^0 (V),$ with tautological bundle $\mathcal E$ of rank $2n$. Consider the product $S^{[n]} \times G (2n, 4n-1)$ with projection
$$\pi: S^{[n]} \times G (2n, 4n-1) \to S^{[n]}.$$
By \cite {F}, Proposition $14.2.2$, $$s_{2n} (V^{[n]}) = \pi_{\star} \left ( c_{\text{top}} (\mathcal E^{\vee} \otimes V^{[n]}) \right).$$ Here $\mathcal E$ and $V^{[n]}$ are pulled back to the product from $G(2n, 4n-1)$ and $S^{[n]}$ respectively. The vector bundle $\mathcal E^{\vee} \otimes V^{[n]}$ has a natural section $s$. 
This is induced by the morphism $$\mathcal E \to V^{[n]} \, \, \text{on} \, \, S^{[n]} \times G (2n, 4n-1)$$ obtained by composing the inclusion $$\mathcal E \hookrightarrow H^0 (V)\otimes \mathcal O$$ on $G(2n, 4n-1)$ with the evaluation map \eqref{ev} $$H^0(V)\otimes \mathcal O\to V^{[n]}\, \, \,\text{on} \, \, \,S^{[n]}.$$  The zero locus $\mathbb D$ of the section $s$ is expected to be zero-dimensional, but let us suppose it is in excess with dimension $d$, so that $$\mathbb D \subset S^{[n]} \times G (2n, 4n-1), \, \, \text{with normal bundle} \, \, N.$$ The excess intersection formula of Section 14.4 in \cite {F} reads 
\begin{equation}
\label{excess0}
s_{2n} (V^{[n]}) = \pi_{\star} \left ( c_{d}  \left (\mathcal E^{\vee} \otimes V^{[n]}|_{\mathbb D}  - N \right )\right), 
\end{equation}
where $\pi: \mathbb D \to S^{[n]}.$

We now compute the right hand side of equation \eqref{excess0}. 
By definition, $$\mathbb D = \{(Z, \Lambda)\, \, {\text{so the map}} \, \, \Lambda \hookrightarrow H^0 (V) \to H^0 (V \otimes \mathcal{O}_Z)\, \, \text{is zero}\}.$$ We describe $\mathbb D$ in concrete terms. Observe that $$V|_{E} = \mathbb C^2 \otimes \mathcal O_{\mathbb P^1}  (n-2),$$ and set $$\Lambda_0 = \text{ker} \, H^0(V) \to H^0 (V|_E). $$ The argument of Lemma \ref{propexpl} shows that $H^1(V\otimes \mathcal O(-E))=0$. Therefore $$\dim \Lambda_0 = (4n-1) - (2n-2) = 2n+1.$$ Observe next that
$$\mathbb D = \{(Z, \Lambda)\, \, {\text{with}} \, \, Z \in E^{[n]} \, \, \text{and} \, \, \Lambda \subset  \Lambda_0\}\subset S^{[n]} \times G(2n, 4n-1).$$ Indeed, if $Z\not \subset E$, by Lemma \ref{propexpl} $$H^0(V)\to H^0(V|_{Z})$$ is surjective. The kernel is of dimension $(4n-1)-2n=2n-1$, hence it cannot contain a subspace $W$ of dimension $2n$. Furthermore, for $Z\subset E$, the restriction $$H^0(V)\to H^0(V|_Z)$$ factors through $H^0(V|_{E})$. Since $V|_E=\mathcal O_E(n-2)\otimes \mathbb C^2$, the restriction $$H^0(V|_{E})\to H^0(V|_{Z})$$ is injective. Therefore, $$W\hookrightarrow H^0(V)\to H^0(V|_{E}) \text{ is zero } \implies \Lambda\subset \Lambda_0.$$
Consequently, we have 
$$\mathbb D = E^{[n]} \times G (2n, \Lambda_0)\subset S^{[n]} \times G (2n, 4n-1),$$ so that $$\mathbb D\simeq \mathbb P^n \times \mathbb P^{2n}.$$ This identification holds scheme-theoretically. Indeed, it is easy to check that the above pointwise arguments can also be carried out in families. The key observation is that cohomology and base change commute for all relative constructions involved.

In $K$-theory, the normal bundle to $\mathbb D$ is 
$$N = N_{E^{[n]}/S^{[n]}}+N_{\mathbb P^{2n}/G(2n, 4n-1)}=\mathcal O (E)^{[n]} + {\mathcal E}^{\vee} \otimes \mathbb C^{2n-2}.$$ Here $2n-2$ is the dimension of the quotient $H^0(V)/\Lambda_0$ and the two summands are restricted to $\mathbb D.$  
We now calculate
\begin{eqnarray*}
s_{2n} (V^{[n]}) &=& \pi_{\star} \left ( c_{3n}  \left ({\mathcal E}^{\vee} \otimes V^{[n]} - \mathcal O (E)^{[n]} - {\mathcal E}^{\vee} \otimes \mathbb C^{2n-2}|_{\mathbb D}   \right )\right)\\
& = & \pi_{\star} \left ( c_{3n}  \left ({\mathcal E}^{\vee} \otimes \mathbb C^2 \otimes \mathcal O(n-2)^{[n]} - \mathcal O (-1)^{[n]} - {\mathcal E}^{\vee} \otimes \mathbb C^{2n-2}   \right )\right) 
\end{eqnarray*}
Let $\zeta$ denote the hyperplane class on $E^{[n]} \simeq \mathbb P^n$ and $h$ denote the hyperplane in the Grassmannian $G(2n, \Lambda_0) \simeq \mathbb P^{2n}.$ We have 
\begin{eqnarray*}
c({\mathcal E}^{\vee}) &=& \frac{1}{1-h},\\
c \left (\mathcal O (n-2)^{[n])} \right ) &=& 1-\zeta,\\
c \left (\mathcal O (-1)^{[n]} \right) &=& (1- \zeta)^n.\\
\end{eqnarray*} The last two formulas for the total Chern classes of tautological vector bundles over the Hilbert scheme of points on $\mathbb P^1$ are explained for instance in the proof of Theorem $2$ of \cite {mop2}.  
Completing the calculation, 
\begin{eqnarray*}
s_{2n} (V^{[n]}) &=& \text{Coeff}_{\,h^{2n} \zeta^n} \left [  (1-\zeta)^{3n+2} \cdot (1-h)^{2n-2} \cdot \frac{1}{c (h^{\vee} \otimes {\mathcal O} (n-2)^{[n]})^2} \right ] \\
&=& \text{Coeff}_{\,h^{2n} \zeta^n} \, \frac{(1-\zeta)^{3n+2} }{(1-h-\zeta)^2} \\
& = & (-1)^n (2n+1).\qed
\end{eqnarray*}

\vskip.2in

\subsection {The case $k=n-1$.} Consider now the case $k=n-1$, so that $$V=V_0\otimes E^{-(n-1)}$$ with $V_0$ satisfying \eqref{apt}, and $\chi(V)=4n-1$. 
By Lemma \ref{propexpl} (ii), all nontrivial extensions \begin{equation}\label{e45}0\to V^{\vee} \otimes E\to W\to I_Z\to 0\end{equation} must have the middle term $W$ a rigid $H$-stable sheaf of rank $3$. Of course, \begin{equation}\label{c1www}c_1(W)=-H+2nE.\end{equation}

We record the following  
\begin{lemma} \label{lw} There are no rigid $H$-stable sheaves of rank $3$ satisfying \eqref{c1www} when $n\not\equiv 0\mod 3$. \vskip.1in
When $n\equiv 0\mod 3$, the sheaf $W$ must necessarily be of the form $$W=\pi^{\star} W_0\otimes E^{\frac{2n}{3}}$$ where $W_0\to X$ is a bundle over the $K3$ surface with
\begin{equation}
\text{rank }W_0=3,\,\,\, c_1 (W_0) = -H, \,\,\, W_0\,\, \,\text{rigid and Gieseker }H\text{-stable}.
\label{apt1}\end{equation}
\end {lemma} 

\proof As noted in \cite {G}, \cite{Y}, stability with respect to the non-ample divisor $H$ on the blowup is equivalent to $H-\epsilon E$-stability for $\epsilon$ sufficiently small. Thus the moduli space of $H$-stable sheaves on the blowup admits a natural scheme structure. In general, the moduli space $M^S_{H}(3, -H+2nE, c_2)$ of $H$-stable sheaves on the blowup has expected dimension $6d-16$ where $d$ denotes the discriminant. For rigid moduli spaces, we have $d=\frac{8}{3}$. The corresponding moduli of $H$-stable sheaves is zero dimensional and consists of isolated points. To count the stable rigid sheaves, we compute the Euler characteristic of the moduli space. 

To this end, we use the blowup formulas of Proposition 3.4 in \cite {Y} or Proposition 3.1 (2) in \cite {G}. The generating series of Euler characteristics on the blowup is computed in terms of the same series on the underlying $K3$ surface $X$ as $$\sum_{d} q^{d-\frac{8}{3}} \mathsf e(M^S_{H}(3, -H+2nE, d)) = \prod_{m=1}^{\infty} \frac{1}{(1-q^m)^{3}}\cdot \left(\sum_{(x, y)} q^{x^2+xy+y^2}\right)$$ $$\cdot \left(\sum_{d} q^{d-\frac{8}{3}} \mathsf e(M_H^X(3, -H, d))\right)$$ where $x, y\in \mathbb Z+\frac{2n}{3}$. To complete the proof, we compute the constant term in the above expression. If $n\not \equiv 0\mod 3$, there is no constant term due to the factor $$\sum q^{x^2+xy+y^2}.$$ Indeed, as $n\not \equiv 0\mod 3$, we have $(x, y)\neq (0, 0)$ so $x^2+xy+y^2>0$. 

For $n\equiv 0\mod 3$, the constant term is $1$. The moduli space consists of a single rigid sheaf. Note that when $\chi (V) = 4n-1,$ a direct calculation shows that there exists a rigid sheaf $W_0$ on $X$ such that $W = \pi^{\star} W_0\otimes E^{\frac{2n}{3}}$ on $S$ has the numerics determined by the exact sequence \eqref{e45}. Specifically the numerical assumptions on $V$ require $$H^2=4n^2+12n-14$$ and we select $W_0$ to be in the moduli space with Chern numbers $$\text{ch}(W_0)=3-H+\left(\frac{2}{3}n^2+2n-5\right)\left[\text{pt}\right].$$ This $W$ is then the only sheaf in its moduli space. 


\qed 
\vskip.1in
By the Lemma, if $n\not \equiv 0\mod 3$ the extension $W$ cannot exist. Consequently, we have $$H^1(V\otimes I_Z)=0$$ for all $Z$. This shows that the evaluation map \eqref{ev} over the Hilbert scheme is surjective, and we have obtained the following result corresponding to Theorem \ref{t4}\,(ii):
\begin{proposition}\label{sub}
If $V=\pi^{\star} V_0\otimes E^{-(n-1)}$ with $\chi(V)=4n-1$ and $V_0$ a rank $2$ bundle as in \eqref{apt},  then $$\int_{S^{[n]}} s_{2n}(V^{[n]})=0 \, \, \, \, \text{for} \, \, \, n\not \equiv 0 \mod 3.$$
\end {proposition} 

We now focus on the remaining statement of Theorem \ref{t4}. We show 
\begin {proposition} \label{p7} If $V=\pi^{\star} V_0\otimes E^{-(n-1)}$ with $\chi(V)=4n-1$ and $V_0$ a rank $2$ bundle as in \eqref{apt},  then $$\int_{S^{[n]}} s_{2n}(V^{[n]})=1 \, \, \, \text{for} \, \, \, n\equiv 0\mod 3.$$ 
\end{proposition} 

\proof  As usual, the top Segre class is supported on the locus $\mathbb D$ where the evaluation map \eqref{ev} is not surjective. We now identify $\mathbb D$ as a subvariety of the Hilbert scheme $S^{[n]}$. The nonsurjectivity locus consists of subschemes $Z$ such that $$H^1(V\otimes I_Z)\neq 0,$$ 
in other words corresponding to the existence of nontrivial extensions
\begin{equation}
\label{basic}
0 \to V^{\vee} \otimes E \to W \to I_Z \to 0. 
\end{equation}
When $n = 3 \ell$ for $\ell \in \mathbb Z,$ as seen in Lemma \ref{lw}, all such extensions have the same middle term, the unique stable rigid rank 3 vector bundle $W$ on $S$ with numerics specified by the lemma. 
A direct calculation establishes further that \begin{equation}\label{chivw}\chi(V^{\vee}\otimes E, W)=1.\end{equation}
We now show that $$\mathbb D \simeq \mathbb P \text{Hom} (V^{\vee} \otimes E, W) = \mathbb P H^0 (V \otimes W \otimes E^{-1}).$$
To start, note the useful vanishing 
$$H^1 (V \otimes W) = 0$$
which will be established in Lemma \ref{l8} using Reider arguments. By Serre duality and stability $$H^2(V\otimes W)=0, \,\, H^2(V\otimes V^{\vee}\otimes E)=\mathbb C.$$ From the exact sequence \eqref{basic}, upon tensoring with $V$, 
these vanishings imply 
\begin{equation}\label{vizc} H^1 (V\otimes I_Z) = \mathbb C.\end{equation} 
By \eqref{vizc}, for any given $Z$ in the degeneracy locus, the basic extension
\eqref{basic} is unique. Conversely, any non-zero morphism $$V^{\vee}\otimes E\to W$$ must be injective by stability, and the cokernel sheaf $Q$ must have rank $1$ and trivial determinant. We claim that $Q=I_Z$ for some scheme $Z$ which necessarily has length $n$. This is clear if $Q$ is torsion-free. However, if 
$Q$ had torsion $T$, then there would be a torsion-free quotient $Q^{\text{tf}}$ of rank $1$, 
$$0 \to T \to Q \to Q^{\text{tf}} \to 0, \, \, \text{and} \, \, Q^{\text{tf}}=I_U(-D), \, \, \text{for a subscheme} \, \,U.$$ Comparing determinants, we find $D=\det T$ effective so that $D\cdot H\geq 0$. This however contradicts stability of $W\to Q^{\text{tf}}$ unless $D\cdot H=0,$ in which case $D=qE$  for some $q\geq 0$. We argue that $q=0$. Indeed, consider the kernel $K$ of the surjection $$W\to Q\to Q^{\text{tf}}=I_U(-qE).$$ Using $K\hookrightarrow W$ and $W$ is $H$-stable, it follows that $K$ is $H$-stable as well. Write $$m=\text{length }(U).$$ A direct calculation shows that \begin{eqnarray*}\chi(K, K)&=&\chi(W-I_U(-qE), W-I_U(-qE))\\&=&2+4m+4n(q-1)+3q^2>2,\end{eqnarray*} which contradicts stability. The only exception is $$q=0,\,\,m\leq n\implies Q^{\text{tf}}=I_U.$$ The torsion part $T$ must have trivial determinant hence it must be supported on points. Computing Euler characteristics $$\chi(I_U)+\chi(T)=\chi(Q)=\chi(V^{\vee}\otimes E)-\chi(W)=1-n\implies \chi(T)=m-n\leq 0.$$ This is only possible if $T=0$ so that $Q$ is torsion free, in fact $Q=I_U$ for some zero dimensional scheme $U$ of length $n$.

Thus we have identified the degeneracy locus $\mathbb D$ as 
$$\mathbb D = \mathbb P H^0 (V \otimes W \otimes E^\vee).$$ It is straightforward to carry out the above argument in families. The crux of the matter is that \eqref{vizc} has constant rank, hence cohomology and base change commute over the degeneracy locus $\mathbb D$.
We set $$ d= \dim \mathbb D.$$ Since $H^2 (V \otimes W \otimes E^\vee) = 0$ by Serre duality and stability, it follows from \eqref{chivw} that $$h^1(V\otimes W\otimes E^{\vee})=d.$$

\vskip.1in

On $\mathbb D \times S$ we have a universal exact sequence
\begin{equation}
\label{basicu}
0 \to V^{\vee} \otimes E \otimes \mathcal O (-1) \to W \to \mathcal I_{\mathcal Z} \to 0. 
\end{equation}
which restricts to \eqref {basic} for each point of $\mathbb D$. 
Here $V, W, E$ are pulled back to the product from $S$, $\mathcal O (1)$ is the hyperplane bundle on $\mathbb D$,
and $\mathcal I_{\mathcal Z}$ is the restriction to $\mathbb D \times S$ of the universal ideal sheaf on $S^{[n]} \times S$.

\vskip.1in

The Segre integral $$\int_{S^{[n]}} s_{2n} (V^{[n]})$$ is calculated as an excess intersection on the degeneracy locus $\mathbb D$ as follows.
Consider the Grassmannian $G (2n, 4n-1)$ of subspaces $$\Lambda\hookrightarrow H^0 (V) = \mathbb C^{4n-1}.$$ Write $\mathcal E$ for the tautological bundle, and let $\mathcal F$ denote the tautological quotient. The vector bundle $$\mathcal E^{\vee}\otimes V^{[n]}\to G (2n, 4n-1)\times S^{[n]}$$ has a natural section $s$ obtained as the composition $$\Lambda\to H^0(V)\to H^0(V\otimes \mathcal O_Z).$$ Let 
$$\mathbb D_0 \subset G (2n, 4n-1) \times S^{[n]}$$ be the zero locus of this section. It consists of those pairs $(Z, \Lambda)$ where $$\Lambda\subset H^0(V\otimes I_Z).$$ In particular, $h^0(V\otimes I_Z)\geq \dim \Lambda = 2n.$
Note now that 
$$\chi(V\otimes I_Z)=\chi(V)-2n=2n-1,$$
while $H^2(V\otimes I_Z)=0$ by Serre duality and stability. 
For $(Z, \Lambda)$ in the degeneracy locus $\mathbb D_0,$ we have then $$h^1(V\otimes I_Z)\geq 1.$$ This shows that $Z$ must be in the degeneracy locus $\mathbb D$, and in this case $h^1(V\otimes I_Z) = 1$ by \eqref{vizc}. We are thus forced to have $\Lambda=H^0(V\otimes I_Z)$ so that  
$$\mathbb D_0 = \{ (Z, \, \Lambda=H^0 (V\otimes  I_Z) ) \, \, \text{with} \, \, Z \in \mathbb D \} \subset S^{[n]} \times G (2n, 4n-1).$$
Note that the projection onto $S^{[n]}$ induces an isomorphism 
$$\mathbb D_0 \simeq \mathbb D.$$

By \cite{F}, Section $14.4$, we have $$\int_{S^{[n]}} s_{2n} (V^{[n]}) = \int_{\mathbb D_0} c_d (\mathbb E),$$ where the excess virtual bundle
$\mathbb E$ is 
$$\mathbb E = \mathcal E^{\vee} \otimes V^{[n]} - N,$$ with $N$ the normal bundle of ${\mathbb D_0} \subset S^{[n]} \times G (2n, 4n-1).$ We have in 
K-theory,
$$N = \mathcal E^{\vee} \otimes \mathcal F + T S^{[n]} - T\mathbb D_0.$$
Putting all together, 
\begin{equation}
\label{excess}
\mathbb E = \mathcal E^{\vee} \otimes (V^{[n]} - \mathcal F) + T\mathbb D_0  - TS^{[n]}.
\end{equation}
It remains to identify $\mathbb E$ explicitly in terms of the universal sequence \eqref{basicu}.
For simplicity we write identities in $K$-theory on fibers. To start, note that on $\mathbb D_0$ we have
\begin{eqnarray*}
\mathcal E|_Z &=& H^0 (V\otimes I_Z), \\
\mathcal F|_Z &= & H^0 (V) - H^0 (V \otimes I_Z),\\ 
 TS^{[n]}|_Z &=& \text{Ext}^1 (I_Z, I_Z),
 \end{eqnarray*}
and from the Euler sequence 
\begin{equation}
 \label{tlocus} 
 T\mathbb D_0 = H^0 (V\otimes W \otimes E^\vee) \otimes \mathcal O (1) - \mathbb C.
\end{equation}

The evaluation sequence 
$$0 \to H^0 (V \otimes I_Z) \to H^0 (V) \to H^0 (V \otimes \mathcal O_Z) \to H^1 (V \otimes I_Z) \to 0$$
gives 
$$V^{[n]} - \mathcal F = R^1\pi_{\star} (V \otimes \mathcal {I_Z}),$$ a line bundle on $\mathbb D_0$. 
Here $\pi: \mathbb D \times S \to \mathbb D$ is the projection. From the sequence \eqref{basicu}, after tensoring with $V$ and pushing forward to $\mathbb D$ via $\pi$ we obtain
$$ 0 \to \mathcal O (-1) \to H^0 (V \otimes W) \to R^0 \pi_{\star}  (V \otimes \mathcal {I_Z}) \to 0,$$
$$ 0 \to R^1 \pi_{\star}  (V\otimes \mathcal {I_Z}) \to \mathcal O (-1) \to 0.$$
Therefore in $K$-theory
\begin{eqnarray*}
R^0 \pi_{\star}  (V \otimes \mathcal {I_Z}) &=& H^0 (V \otimes W) - \mathcal O (-1) \\
R^1 \pi_{\star} (V \otimes \mathcal {I_Z}) &=& \mathcal O (-1).
\end{eqnarray*}
We obtain
\begin{equation}
\mathcal E^{\vee} \otimes (V^{[n]} - \mathcal F) =  H^0 (V \otimes W)^{\vee} \otimes \mathcal O (-1) - \mathbb C. \label{grass} 
\end{equation}

Noting from \eqref{basicu} that $$\mathcal I_{\mathcal Z} = W - V^{\vee} \otimes E \otimes \mathcal O (-1),$$ 
it follows that $$\mathcal I_{\mathcal Z}^{\vee} \, \otimes \, \mathcal I_{\mathcal Z} = V \otimes V^{\vee} + W \otimes W^{\vee} - V \otimes W \otimes E^\vee \otimes \mathcal O (1) - V^{\vee}  \otimes W^{\vee} \otimes E \otimes \mathcal O (-1).$$
Here all tensor products are derived.
We therefore finally
calculate
\begin{eqnarray}
\label{tangent}
TS^{[n]} &=&\, {\mathcal{E}xt}^1_{\pi} (\mathcal {I_Z}, \, \mathcal {I_Z}) \\ \nonumber &=&-{\mathcal{E}xt}^{\bullet}_{\pi} (\mathcal {I_Z},  \mathcal {I_Z})+\mathbb C+\mathbb C\\ \nonumber&=&  H^{\bullet} (V \otimes W \otimes E^\vee) \otimes \mathcal O (1)  + H^0 (V \otimes W)^{\vee} \otimes \mathcal  O (-1) -   \mathbb C - \mathbb C. 
\end{eqnarray} This last equality uses the fact that $V, W$ are stable and rigid, and that $V\otimes W$ has no higher cohomology. 

Collecting the expressions \eqref{tlocus}, \eqref{grass}, \eqref{tangent} in the excess bundle $\mathbb E$ given by \eqref{excess},
we find 
$$\mathbb E = H^1 (V \otimes W \otimes E^\vee) \otimes \mathcal O (1) = \mathbb C^d \otimes \mathcal O (1).$$
We conclude 
$$\int_{X^{[n]}} s_{2n} (V^{[n]}) = \int_{\mathbb D_0} c_d (\mathbb E) =  \int_{\mathbb D_0} h^d = 1.$$
This completes the proof. 
\qed

\begin{lemma}
\label{l8} 
Let $n=3\ell$. Assume $$V=\pi^{\star}V_0\otimes E^{-3\ell+1}, \,\,W=\pi^{\star}W_0\otimes E^{2\ell}$$ 
where $\chi(V)=4n-1$ and $V_0, W_0$ are two rigid bundles satisfing conditions \eqref{apt} and \eqref{apt1}. We have $$H^1(V\otimes W)=0.$$
\end{lemma} 

\proof The argument is an application of Reider's method. Assume that $$H^1(V\otimes W)=\text{Ext}^{1}(W, V^{\vee}\otimes E)^{\vee}\neq 0.$$ We construct a nontrivial extension $$0\to V^{\vee}\otimes E\to F\to W\to 0.$$ The middle sheaf $F$ is $H$-stable by an argument similar to that of Proposition \ref{p3}. 
Indeed, if $F$ is not $H$-stable, let $G \hookrightarrow F$ be the maximal semistable destabilizing subsheaf. Then $${\text{rk}\, G} < {\text{rk}\, F} \, \, \, \text{and} \, \,  \,   \mu_H(G) \geq \mu_H(F)>\mu_H(V^{\vee}\otimes E)\,.$$
We see that $G$ cannot be a subsheaf of the kernel $V^{\vee}\otimes E$ since this would contradict the $H$-stability of the latter. Therefore, we have a nonzero morphism $$\phi: G \to W\implies \mu_H(G)<\mu_H(W).$$ Writing $c_1(G)=aH+bE$ we see that $$\mu_H(F)\leq\mu_H(G)<\mu_H(W)\implies -\frac{2}{5}\leq \frac{a}{\text{rk}\,G}<-\frac{1}{3}$$ which is impossible. 

Finally, a direct calculation shows \begin{eqnarray*}\chi(F, F)&=&\chi(W, W)+\chi(V^{\vee}\otimes E, V^{\vee}\otimes E)+\chi(V^{\vee}\otimes E, W)+\chi(W, V^{\vee}\otimes E)\\&=&2+2+1+(2n+1)>2\end{eqnarray*} contradicting the  stability of $F$. \qed

\subsection {Proof of Theorem \ref{conj2}.} We now prove Theorem \ref{conj2}. The statement will follow combinatorially using the geometric input provided by Theorem \ref{t4}.
\vskip.1in

Define first the following combination of the basic power series $A_i$: \begin{equation}\label{define}f(z)=A_0(z)^5\cdot A_1(z)^{20}\cdot A_3(z)^2\end{equation}  
$$g(z)=A_0(z)^{-4} \cdot A_1(z)^{-22}\cdot A_2(z)^2\cdot A_3(z)^{-4} \cdot A_4(z)^{-1}$$ $$h(z)=A_0(z)^{-3}\cdot A_1(z)^{-18}\cdot A_2(z)^2\cdot A_3(z)^{-2}\cdot A_4(z)^{-1}.$$ We will derive identities between the functions $f, g, h$ using the two calculations provided by Theorem \ref{t4}.

First, over the $K3$ blowup, let $V=V_0\otimes E^{-(n-2)}$ with $V_0$ rigid $H$-stable of rank $2$ so that $\chi(V)=4n-1$. \footnote{As noted in Lemma \ref{propexpl}, this numerical setup only exists over $K3$s of certain genus determined by $n$. The argument here strongly uses the universality of the Segre series.} One checks that $$c_1(V)^2=20n-22,\,\, c_2(V)=5n-4, \,\, c_1(V)\cdot K_S=2(n-2), \,\,K_S^2=-1.$$
Using \eqref{formm2} for the vector bundle $V$, we obtain that $$\int_{S^{[n]}} s(V^{[n]})=\left[z^n\right] A_0(z)^{5n-4} \cdot A_1(z)^{20n-22} \cdot A_2(z)^2\cdot A_3(z)^{2n-4}\cdot A_4(z)^{-1}=f(z)^n \cdot g(z).$$ With the aid of Proposition \ref{asp} this rewrites as 
$$\left[z^n\right] f(z)^{n}\cdot g(z)=(-1)^n (2n+1),$$ where the brackets denote the suitable coefficient in the given power series. Therefore $$\sum_{n=0}^{\infty} \left[z^n\right] f(z)^{n}\cdot g(z)=\sum_{n=0}^{\infty} z^n (-1)^n(2n+1)=\frac{1-z}{(1+z)^2}.$$ 

Write \begin{equation}\label{fw}z=\frac{w}{f(w)}.\end{equation} The Lagrange-B\"urmann inversion formula \cite {WW} is the general identity $$\sum_{n=0}^{\infty} \left(\left[z^n\right] f(z)^{n}\cdot g(z)\right)\cdot z^n=\frac{g(w)}{f(w)} \cdot \frac{dw}{dz}.$$ In our situation, it gives \begin{equation}\label{first1}\frac{1-z}{(1+z)^2}=\frac{g(w)}{f(w)} \cdot \frac{dw}{dz}.\end{equation}

In similar fashion, for $V=V_0\otimes E^{-(n-1)}$, by making use of Proposition \ref{sub} and Proposition \ref{p7}, we obtain \begin{equation}\label{2nd}\frac{1}{1-z^3}=\frac{h(w)}{f(w)} \cdot \frac{dw}{dz}.\end{equation} The expression $$\frac{1}{1-z^3}=\sum_{k=0}^{\infty} z^{3k}$$ encodes the fact that the Segre integrals are $0$ for $n\not \equiv 0\mod 3$ and equal to $1$ for $n\equiv 0 \mod 3$. 

We now explain how equations \eqref{first1} and \eqref{2nd} give the remaining functions $A_3$ and $A_4$. 
Dividing the two equations we obtain $$\frac{h(w)}{g(w)}=\frac{(1+z)^2}{(1-z)(1-z^3)}.$$ This gives via \eqref{fw} $$\frac{w}{f(w)} \cdot \frac{h(w)}{g(w)}=z\cdot \frac{(1+z)^2}{(1-z)(1-z^3)}.$$ Let us write $$w=t(1+3t)^3.$$ We compute 
$$\frac{w}{f(w)} \cdot \frac{h(w)}{g(w)}=w\cdot A_0(w)^{-4}\cdot A_1(w)^{-16}=\frac{t}{1+3t},$$ where the first equality follows by \eqref{define}, and the second equality uses the expressions for $A_0, A_1$ given by Theorem \ref{t1}.  
Therefore $$\frac{z(1+z)^2}{(1-z)(1-z^3)}=\frac{t}{1+3t}\implies z=\mathsf y(t),$$ for the function $\mathsf y(t)$ of the introduction. With this understood, we find via \eqref{fw} $$f(w)=\frac{w}{z}=\frac{t(1+3t)^3}{\mathsf y}.$$ Using \eqref{define}, we obtain $$A_3(w)= f(w)^{\frac{1}{2}} \cdot A_0(w)^{-\frac{5}{2}}\cdot A_1(w)^{-10}=\frac{1}{1+3t}\cdot \left(\frac{t}{\mathsf y}\right)^{1/2}$$ where in the last equality we used the expressions for $A_0, A_1$ in Theorem \ref{t1}.  

Similarly, from equation \eqref{first1} we compute $$g(w)=f(w)\cdot \frac{1-z}{(1+z)^2}\cdot \frac{dz/dt}{dw/dt}=\frac{t(1+3t)^3}{\mathsf y}\cdot \frac{1-\mathsf y}{(1+\mathsf y)^2}\cdot \frac{\mathsf y'}{(1+3t)^2(1+12t)}.$$ Combined with \eqref{define} we obtain the expression for $A_4$ claimed in Theorem \ref{conj2}. \qed

\end{document}